\documentclass[11pt]{article}
\usepackage{amstex,amssymb,latexsym}
\usepackage{arthead}
\usepackage{mpmath}
\newtheorem{Theorem and Definition}[Definition]{Theorem and Definition}
\newtheorem{Proposition and Definition}[Definition]{Proposition and Definition}

\newenvironment{rExample}{\begin{Example} \rm}{ \end{Example}}

\newenvironment{rNote}{\begin{Note} \rm}{ \end{Note}}
\newcommand{\pse}[1]{\Psi_{\rho,\delta}^{#1}}
\newcommand{\psei}{\Psi_{\rho,\delta}^{\infty}}
\newcommand{\psemi}{\Psi^{-\infty}}
\renewcommand{\sym}[1]{{\rm S}_{\rho,\delta}^{#1}}
\renewcommand{\symi}{{\rm S}_{\rho,\delta}^{\infty}}

\renewcommand{\symmi}{{\rm S}^{-\infty}}

\newcommand{\ddera}[1]{\frac{D^{|\alpha|}}{\partial {#1}^{\alpha}}}
\newcommand{\dderb}[1]{\frac{D^{|\beta|}}{\partial {#1}^{\beta}}} 

\newcommand{\ddere}[2]{\frac{D^{|{#1}|}}{\partial {#2}^{#1}}}
\begin{document}
\thispagestyle{empty}
\title{The normal symbol on Riemannian manifolds}
\author{ Markus J.~Pflaum\thanks{ Humboldt Universit\"at zu Berlin, 
         Mathematisches Institut, Unter den Linden 6, 10099 Berlin}}
\date{ \today }
\maketitle
\begin{abstract}
For an arbitrary Riemannian manifold $X$ and Hermitian vector bundles $E$ and
$F$ over $X$ we define the notion of the normal symbol of a 
pseudodifferential operator $P$ from $E$ to $F$.
The normal symbol of $P$ is a certain smooth function from the cotangent 
bundle $T^*X$ to the homomorphism bundle ${\rm Hom} (E,F)$
and depends on the metric structures resp.~the corresponding connections 
on $X$, $E$ and $F$.   
It is shown that by a natural integral formula the 
pseudodifferential operator $P$ can be recovered from  its symbol.
Thus, modulo smoothing operators resp.~smoothing symbols, we receive a linear 
bijective correspondence between the space of symbols and the space 
of pseudodifferential operators on $X$. This correspondence
comprises a natural transformation between appropriate functors.
A formula for the asymptotic expansion of the product symbol 
of two pseudodifferential operators in terms of the symbols of its
factors is given. Furthermore an expression for the symbol of the
adjoint is derived.
Finally the question of invertibility of pseudodifferential operators
is considered. For that we use the normal symbol to establish a new and 
general notion of elliptic pseudodifferential operators on manifolds.
\end{abstract}

\tableofcontents
\mbox{ } 
{\bf keywords:} pseudodifferential operators on manifolds, 
                asymptotic expansions, symbol calculus


%
%
%
%
\section*{Introduction}
\addcontentsline{toc}{section}{Introduction}
It is a well-known fact that on Euclidean space one can construct a
canonical linear isomorphisms between symbol spaces and corresponding
spaces of pseudodifferential operators. Furthermore one has natural
formulas which represent a pseudodifferential operator in terms of
its symbol resp.~which give an expression for the symbol of a
pseudodifferential operator. By using symbols one gets much
insight in the structure of pseudodifferential operators on Euclidean
space; in particular they give the means to construct a 
{(pseudo)} inverse of an elliptic differential operator on $\R^d$.

Compared to the symbol calculus for pseudodifferential operators on $\R^n$ 
it seems that the symbol calculus for pseudodifferential operators on manifolds
is not that well-established.
But as the pure consideration of symbols and operators in local
coordinates does not reveal the geometry and topology of the manifold 
one is working on, it is very desirable to build up a general theory of 
symbols for pseudodifferential operators on manifolds.
In his articles \cite{Wid:FPO,Wid:CSCPO}  {\sc Widom}
gave a proposal for a symbol calculus on manifolds.
By using a rather general notion of a phasefunction {\sc Widom} constructs
a map from the space of pseudodifferential operators on manifolds to
the space of symbols and shows by an abstract argument for the case
of scalar symbols that this map is bijective modulo smoothing operators
resp.~symbols.

In our paper we introduce a symbol calculus for pseudodifferential operators
between vector bundles having the feature that both the symbol map and
its inverse have a concrete representation. 
In particular we thus succeed in giving an integral representation for the 
inverse of our symbol map or in other words for the operator map. 
Essential for our approach is an appropriate notion of a phasefunction.
Because of our special choice of a phasefunction the resulting symbol calculus 
is natural in a category theoretical sense.

Using the integral representation for the operator map it is possible to 
write down a formula for the symbol of the adjoint of a pseudodifferential 
operator and for the symbol of the product of two pseudodifferential operators.

The normal symbol calculus will furthermore give us the means to 
build up a natural notion of elliptic symbols respectively elliptic 
pseudodifferential operators on manifolds. 
It generalizes the classical notion of ellipticity as defined for example in 
{\sc H\"ormander} \cite{Hor:ALPDOIII} and also the concept of ellipticity 
introduced by {\sc Douglis, Nirenberg} \cite{DouNir:IEESPDE}. 
Our framework of ellipticity allows the construction of parametrices of 
nonclassical respectively nonhomogeneous elliptic pseudodifferential
operators. Moreover we do not need principal symbols for defining ellipticity.
Instead we define elliptic operators in terms of their normal symbol, 
as the globally defined normal symbol of a pseudodifferential operator 
carries more information than its principal symbol.
Moreover the normal symbol of a pseudodifferential operator 
shows immediately, whether the corresponding operator is invertible 
modulo smoothing operators or not.     

Let us also mention that another application of the normal symbol calculus 
on manifolds lies in quantization theory. 
There one is interested in a quantization map associating pseudodifferential
operators to certain functions on a cotangent bundle in a way that
{\sc Diracs} quantization condition is fulfilled. But the inverse of a 
symbol map does exactly this, so it is a quantization map.
See {\sc Pflaum} \cite{Pfl:LADQ} for details.

Note that our work is also related to the recent papers of 
{\sc Yu.~Safarov} \cite{Saf:POLC,Saf:FLBO}.


\section{Symbols}
First we will define the notion of a symbol on a Riemannian vector bundle.
Let us assume once and for all in this article that $\mu, \rho, \delta \in \R$ 
are real numbers such that $0 \leq \delta < \rho \leq 1$ and 
$1 \leq \rho + \delta$. The same shall hold for triples
$\tilde{\mu}, \tilde{\rho},\tilde{\delta} \in \R$.

We sometimes use properties of symbols on $\R^d \times \R^N$. As these are
well-described in the mathematics literature we only refer the interested
reader to  {\sc H\"ormander} \cite{Hor:ALPDOI,Hor:ALPDOIII}, 
{\sc Shubin} \cite{Shu:POST} or 
{\sc Grigis, Sj\"ostrand} \cite{GriSjo:MADO}  for a general introduction to 
symbol theory on $\R^d \times \R^N$ and for proofs.
\begin{Definition}
Let $ E \rightarrow X$ be a Riemannian or Hermitian vector bundle
over the smooth manifold $X$, $\varrho: TX \rightarrow X$ its tangent 
bundle and $\pi : T^*X \rightarrow X$ its cotangent bundle.
Then an element $a \in {\cal C}^{\infty} (\pi^* (E |_U))$
with $U \subset X$ open is called 
a {\bf symbol} \label{SySpVeBu} over $U \subset X$ with {\bf values} in $E$
of {\bf order} $\mu$ and {\bf type} $(\rho,\delta)$ if for every trivialization
$(x,\zeta) : T^*X|_V \rightarrow \R^{2d} $ and every triviliazation 
$\Psi : E|V \rightarrow V \times \R^N$ with $V \subset U$ open
the following condition is satisfied:
\vspace{0.2cm} \\
(Sy)
For every $\alpha, \beta \in \N^d$ and every $K \subset V$ compact
there exists $C = C_{K,\alpha,\beta} > 0 $ such that for every 
$\xi \in T^*X|_K$ the inequality
\begin{equation}
  \left| \left| \dera{x} \, \derb{\zeta} \, \Psi (a (\xi)) \right| \right| 
  \: \leq \: C \, (1 + ||\xi||)^{
  \mu + \delta |\alpha| - \rho |\beta|}
\end{equation}
is valid.
\\[2mm]
The space of these symbols is denoted by 
${\rm S}^{\mu}_{\rho,\delta}(U,T^*X,E)$
or shortly ${\rm S}^{\mu}_{\rho,\delta} (U, E)$. It gives rise to the sheaf
$\sym{\mu} (\cdot,T^*X, E)$ of symbols on $X$ with values in $E$ 
of order $\mu$ and type $(\rho,\delta)$.

By defining ${\rm S}^{-\infty} (U,E) = \bigcap_{\mu \in \R} \,
{\rm S}^{\mu}_{\rho,\delta} (U,E) $ and 
${\rm S}^{\infty}_{\rho,\delta} (U,E) = \bigcup_{\mu \in \R} \,
{\rm S}^{\mu}_{\rho,\delta} (U,E) $ we get 
the sheaf ${\rm S}^{-\infty} (\cdot,E)$ of smoothing symbols
resp.~the sheaf ${\rm S}^{\infty}_{\rho,\delta} (\cdot,E)$ of symbols of type
$(\rho,\delta)$.

In case $E$ is the trivial bundle $X \times \C$ of a Riemannian manifold $X$
we write  $\sym{\mu} $, $\symi $ and $\symmi $ for the corresponding 
symbol sheaves $\sym{\mu} (\cdot,X\times \C)$, $\symi (\cdot, X \times C)$ and 
$\symmi (\cdot,X \times \C)$.

The space of symbols $\sym{\mu} (U,TX,E)$ consists of smooth functions 
$a \in {\cal C}^{\infty} (\varrho^*(E|_U))$ such that $a$ fulfills 
condition (Sy). It gives rise to sheaves $\sym{\mu} (\cdot,TX,E)$, 
${\rm S}^{\infty}_{\rho,\delta} (\cdot,TX,E)$ 
and ${\rm S}^{-\infty}(\cdot,TX,E)$.
\end{Definition}
It is possible to extend this definition to one of symbols over conic manifolds
but we do not need a definition in this generality and refer the 
reader to {\sc Duistermaat} \cite{Dui:FIO} and 
{\sc Pflaum} \cite{Pfl:LADQ}. \vspace{0.2cm}

We want to give the symbol spaces  ${\rm S}^{\mu}_{\rho,\delta} (U, E)$
a topological structure. So first choose a compact set $K \subset U$, a 
(not necessarily disjunct) partition $ K = \bigcup \limits_{\iota \in J} 
K_{\iota}$ of $K$ into compact subsets together with local trivializations
$(x_{\iota} , \zeta_{\iota} ) : T^*X|_{V_{\iota}} \rightarrow \R^{2d}$
and trivializations $\Psi_{\iota} : E|{V_{\iota}} \rightarrow V_{iota} \times 
\R^N$ such that $K_{\iota} \subset V_{\iota} \subset U$ and $V_{\iota}$ open.
Then we can attach for any $\alpha, \beta \in \R^d$ a seminorm
$p=p_{(K_{\iota},(x_{\iota},\zeta_{\iota})),\alpha,\beta} : 
{\rm S}^{\mu}_{\rho,\delta} (U,E) \rightarrow \R^+ \cup \{ 0 \}$ 
to the symbol spaces $ {\rm S}^{\mu}_{\rho,\delta} (U, E) $ by
\begin{equation}
  p (a) \: = \: \sup \limits_{\iota \in J} \left\{
  \frac{\left| \left|  \dera{x_{\iota}} \derb{\zeta_{\iota}} \Psi_{\iota}
  (a (\xi)) 
  \right| \right|}{(1 + |\xi|)^{\mu + |\alpha| \delta - |\beta| \rho} } 
  \, : \: \xi \in T^*X|_{K_{\iota}} \right\}.
\end{equation}
The system of these seminorms gives $ {\rm S}^{\mu}_{\rho,\delta} (U, E) $
the structure of a Fr\'echet space such that the restriction morphisms 
$ {\rm S}^{\mu}_{\rho,\delta} (U, E) \rightarrow
{\rm S}^{\mu}_{\rho,\delta } (V,E) $ for $V \subset U$ open are continuous.
Additionally we have natural and continuous inclusions 
$ {\rm S}^{\mu}_{\rho,\delta} (U, E) \subset 
{\rm S}^{\tilde{\mu}}_{\tilde{\rho},\tilde{\delta}} (U, E) $
for $\tilde{\mu} \geq \mu$, $\tilde{\rho} \leq \rho$ and 
$\tilde{\delta} \geq \delta$. Like in the case of symbols on $\R^d \times \R^N$ 
one can show that pointwise multiplication of symbols is continuous with
respect to the Fr\'echet topology on $ {\rm S}^{\mu}_{\rho,\delta} (U, E) $.
\begin{Proposition}
  Let $a \in {\cal C}^{\infty} (\pi^*(E|_U))$. 
  If there exists an open covering of
  $U$ by patches $V_{\iota}$ with local trivializations $(x_{\iota}, 
  \zeta_{\iota}) : T^*X|_{V_{\iota}} \rightarrow \R^{2d}$
  and trivializations $\Psi_{\iota} : E|{V_{\iota}} \rightarrow V_{\iota} 
  \times \R^N$
  such that for every $(x_{\iota},\zeta_{\iota})$ the above condition (Sy)
  holds, then $a$ is a symbol of order $\mu$ and type $(\rho,\delta)$.
\end{Proposition}
\begin{Proof}
  Let $(x, \zeta) : T^*X|_V \rightarrow \R^{2d}$,
  $V \subset U$ open be a trivialization. Then we can write on $V_{\iota} \cap 
  V$
\begin{equation}
  \dera{x} \derb{\zeta} \: = \: \sum \limits_{\tilde{\alpha} + \alpha' \leq 
  \alpha \atop \iota } 
  \, \left(f_{\tilde{\alpha}, \alpha',\beta} \circ \pi \right) 
  \, \zeta^{\alpha'}_{\iota} \, \dere{\tilde{\alpha}}{x_{\iota}}
  \dere{\beta + \alpha'}{\zeta_{\iota}}\:,
\end{equation}
  where $f_{\tilde{\alpha}, \alpha',\beta} \in {\cal C}^{\infty} (V_{\iota}
  \cap V)$. As $\rho + \delta \geq 1$ holds and (Sy) is true for 
  $(x_{\iota},\zeta_{\iota})$, the claim now follows.
\end{Proof}

\begin{rExample}
\label{ExSy}
\begin{enumerate}  
\item
  Let $X$ be a Riemannian manifold. Then every smooth function 
  $a :T^*U \rightarrow \C$ which is a polynomial function  of order $\leq \mu$
  on every fiber lies in ${\rm S}^{\mu}_{1,0} (U)$.
  In particular we have an embedding ${\cal D}_0 \rightarrow \symi$, 
  where ${\cal D}_0$ is the sheaf of polynomial symbols on 
  $X$, i.~e.~for every $U \subset X$ open ${\cal D}_0 (U) = 
  \{ f: T^*U \rightarrow \C: \: f|_{T^*_xX} \mbox{ is a polynomial for every } 
  x \in U \}$.
\item
  Again let $X$ be Riemannian. Then the mapping 
  $l : T^*X \rightarrow \C$, $\xi \mapsto || \xi ||^2 $ is a symbol of order 
  $2$ and type $(1,0)$ on $X$, but not one of order $\mu < 2$.
  Next regard the function $a : T^*X \rightarrow \C$, 
  $\xi \mapsto \frac{1}{1 + || \xi ||^2} $. This function
  is a symbol of order $-2$ and type $(1,0)$ on $X$, but not one of order
  $\mu < -2$.
\item
  Assume $\varphi: X \rightarrow \R$ to be smooth and bounded. Then 
  $a_{\varphi} :  T^*X
  \rightarrow \R$, $\xi \mapsto (1 + || \xi ||^2 )^{\varphi (\pi (\xi))}$
  comprises a symbol of order $\mu = \sup_{x \in X} \varphi (x)$ and
  type $(1, \delta)$, where $0 < \delta < 1$ is arbitrary.
\end{enumerate}
\end{rExample}

The following theorem is an essential tool for the use of symbols in the
theory of partial differential equations and extends a well-known result for
the case of $E = \R^d \times \R^N$ to arbitrary Riemannian vector bundles.
\begin{Theorem}
\label{MaThAsSuSy}
Let $a_j \in \sym{\mu_j} (U,E)$, $j \in \N$ be symbols such that 
$(\mu_j)_{j \in \N}$ is a decreasing sequence with 
$\lim \limits_{j \rightarrow \infty} \mu_j = - \infty$.
Then there exists a symbol $a \in \sym{\mu_0} (U,E)$ unique up to
smoothing symbols, such that
$ a - \sum \limits_{j=0}^{k} a_j \in \sym{\mu_k} (U,E) $ for every $k \in \N$.

This induces a locally convex Hausdorff topology on the vector space
${\rm S}^{\infty}_{\rho,\delta}/{\rm S}^{-\infty}  (U, E)$,
which is called the {\bf topology of asymptotic convergence}.
\end{Theorem}
\begin{Proof}
  It is a well-known fact that the claim holds for $U \subset \R^d$ and\
  $E = \R^d \times \R^N$ (see \cite{Hor:ALPDOIII,Shu:POST,GriSjo:MADO}). 
  Covering $U$ by trivializations $(z_{\iota},\zeta_{\iota}): T^*X|_{V_{\iota}}
  \rightarrow \R^{2d}$ and 
  $\Psi_{\iota} : E|{V_{\iota}} \rightarrow V_{\iota} \times \R^N$ 
  we can find a partition $(\varphi_{\iota})$
  of unity subordinate to the covering 
  $(V_{\iota})$ and for every index $\iota$ a symbol
  $b_{\iota,k} \in \sym{\mu} (V_{\iota}, \R^N)$ such that 
  $b_{\iota,k} = \sum \limits_{j=0}^k \Psi_{\iota} \circ a_j |_{V_{\iota}} \in 
  \sym{\mu_k} (V_{\iota}, \R^N )$. Now define 
  $a = \sum \limits_{\iota} (\varphi_{\iota} \circ \pi) \, \Psi^{-1}_{\iota}
  \circ b _{\iota,k}$
  and check that this $a$ satisfies the claim. Uniqueness of $a$ up to
  smoothing symbols is clear again from a local consideration.
\end{Proof}
Polynomial symbols are not affected by smoothing symbols. 
The following proposition gives the precise statement.
\begin{Proposition}
\label{EmPoObSySp}
  Let $X$ be a smooth manifold, and ${\cal D}_0 $ the sheaf of polynomial 
  symbols on $X$. Then by composition with the projection
  $\symi \rightarrow \symi / \symmi $ the canonical embedding
\begin{eqnarray}
\begin{array}{cccc}
   {\cal D}_0 & \rightarrow & \symi  \\
   {\cal D}_0 (U) \ni f & \mapsto & f \in \symi (U), & \quad
   U \subset X \mbox{ open}
\end{array}
\end{eqnarray}
  gives rise to a monomorphism 
  $\delta: \: {\cal D}_0 \rightarrow \symi / \symmi$ of sheaves of algebras.
\end{Proposition}
\begin{Proof}
  If $f,g \in {\cal D}_0(U)$ are polynomial symbols such that
  $f -g \in \symmi(U)$, then $f-g$ is a bounded polynomial
  function on each fiber of $T^*U$. Therefore $f-g$ is constant
  on each fiber, hence $f-g = 0$ by $f -g \in \symmi(U)$.
\end{Proof}


\section{Pseudodifferential operators on manifolds}
\label{FoTrPsDiOpMa}
Let $a$ be a polynomial symbol defined on the (trivial) 
cotangent bundle $T^*U = U \times \R^d $ of an open set 
$U \subset \R^d$ of Euclidean space.
Then one can write
$a \: = \: \sum_{\alpha} \, (a_{\alpha} \circ \pi) \zeta^{\alpha}$, 
where $\zeta: T^*U \rightarrow \R^d$ is the projection on the 
``cotangent vectors'' and the $a_{\alpha}$ are smooth functions on $U$. 
According to standard results of partial differential equations one knows 
that the symbol $a$ defines the differential operator $A=\Op{\,(a)}$:  
\begin{equation}
  {\cal C}^{\infty} (U) \ni f \mapsto \sum_{\alpha} 
  \, a_{\alpha} \left(-i \right)^{| \alpha |} 
  \dera{z} f \in {\cal C}^{\infty} (U).
\end{equation}
Sometimes $A$ is called the quantization of $a$. In case $f \in {\cal D} (U)$,
the Fourier transform of $f$ is well-defined and provides 
the following integral representation for $Af$:
\begin{equation}
\label{InRePoOb}
  Af \: = \: \Op{\, (a)} \, (f)  \: = \: 
  \int_{\R^d} \, e^{i\inprod{ \xi}{\cdot}} \, a({\xi})\,\hat{f} (\xi) \, d\xi.
\end{equation}
It would be very helpful for structural and calculational considerations  
to extend this formula to Riemannian manifolds.
To achieve this it is necessary to have an appropriate notion of 
Fourier transform on manifolds. In the following we are going to define such a 
Fourier transform and will later get back to the problem of an 
integral representation for the ``quantization map'' ${\rm Op}$ on Riemannian 
manifolds. \vspace{0.5cm} 
%

Assume $X$ to be Riemannian of dimension $d$ 
and consider the exponential function $\exp$ 
with respect to the  Levi-Civita connection on $X$. Furthermore let
$E \rightarrow X$ be a Riemannian or Hermitian vector bundle.
Choose an open neighborhood $W \subset TX$ 
of the zero section in $TX$ such that 
$(\varrho,\exp{}):W \rightarrow X \times X $
maps $W$ diffeomorphically onto an open neighborhood of the diagonal 
$\Delta$ of $X \times X$. Then there exists a smooth function
$\psi : TX \rightarrow [ 0,1]$ called a {\bf cut-off function}
such that $ \psi |_{\tilde{W}} = 1$
and $\supp{\psi} \subset W$ for an open neighborhood $\tilde{W} \subset
W$ of the zero section in $TX$.  
Next let us consider the unique torsionfree metric connection on $E$.
It defines a parallel transport $\tau_{\gamma} : E_{\gamma (0)} \rightarrow
E_{\gamma_(1)}$ for every smooth path $\gamma: [0,1] \rightarrow X$.
For every $v \in TX$ denote by $\tau_{\exp v}$ or just $\tau_{y,x}$ with
$y = \exp v $ and $x = \varrho (v)$ the parallel transport
along $[0,1] \ni t \mapsto \exp tv \in X$.
Now the following microlocal lift ${\cal M}_{\psi}$ is well-defined:
\begin{eqnarray}
   {\cal M}_{\psi} : & {\cal C}^{\infty} (U, E) \rightarrow \symmi (U,TX,E), &
   f \mapsto \li{f}{\psi}, \\
   & & {\li{f}{\psi}} \, (v) = \left\{ 
\begin{array}{ll}
   \psi (v) \, \tau_{\exp v}^{-1} (f (\exp{v} )) & \mbox{ for } v \in W \\
   0                      & \mbox{ else}.
\end{array} 
\right. 
\end{eqnarray}
${\cal M}_{\psi}$ is a linear but not multiplicative map between function 
spaces.

Over the tangent bundle $TX$ we can define the {\bf Fourier transform} 
as the following sheaf morphism:
\begin{equation}
   {\cal F} : \symmi (U, TX, E ) \rightarrow \symmi (U, T^*X, E) , 
   \quad
   a \mapsto \hat{a} (\xi) = 
   \ipif \int_{T_{\pi (\xi)} X} \! e^{-i \inprod{\xi}{ v}} \, a(v) \, dv.
\end{equation}
We also have a {\bf reverse Fourier transform}:
\begin{equation}
   {\cal F}^{-1} : \symmi (U, T^*X , E) \rightarrow \symmi (U, TX, E) , \quad
   b \mapsto \check{b} (v) = 
   \ipif \int_{T^*_{\pi (v)} X} \! e^{i \inprod{\xi}{v}} \, b (\xi) \, d\xi.
\end{equation}
It is easy to check that
${\cal F}^{-1}$ and ${\cal F} $ are well-defined and inverse to each other 
indeed.

Composing ${\cal M}_{\psi}$ and ${\cal F}$ gives rise to the 
{\bf Fourier transform} ${\cal F}_{\psi} = {\cal F} \circ {\cal M}_{\psi}$
on the Riemannian manifold $X$. ${\cal F}_{\psi}$ has the left inverse
\begin{equation}
  \symmi (U, T^*X,E) \rightarrow {\cal C}^{\infty} (U,E),
  \quad  a \mapsto {\cal F}^{-1} (a) |_{0_U},
\end{equation} 
where $|_{0_U}$ means the restriction to the natural embedding of
$U$ into $T^*X$ as zero section.

By definition ${\cal M}_{\psi}$ and ${\cal F}_{\psi}$ do depend 
on the smooth cut-off function $\psi$, but this arbitrariness
will only have minor effects on the study of pseudodifferential operators
defined through ${\cal F}_{\psi}$.
We will get back to this point later on in this section. \vspace{0.5cm}

The exponential function $\exp$ gives rise to normal coordinates 
$z_x : V_x \rightarrow \R^d$, where $x \in V_x \subset X$, $V_x$ open and 
\begin{equation}
\label{NoCoOrFr}
  \exp_x^{-1} (y) \: = \: \sum \limits_{k=1}^{d} \, z_x^k (y) \cdot e_k (y)
\end{equation}
for all $y \in V_x$ and an orthonormal frame $(e_1,...,e_d)$ of $TV_x$.
Furthermore we receive bundle coordinates $(z_x,\zeta_x): T^*V_x \rightarrow 
\R^{2d}$ and $(z_x, v_x): TV_x \rightarrow \R^{2d}$.
Note that for fixed $x \in X$ the map
\begin{equation}
  V_x \times T^*V_x \ni (y , \xi) \mapsto \left( z_y(\pi (\xi)) , \zeta_y (\xi) 
  \right) \in \R^{2d}
\end{equation}
is smooth and that by the Gau{\ss} Lemma 
\begin{equation}
\label{SyGaLe}
  z_x^k(y) = - z_y^k(x)
\end{equation}
for every $y \in V_x$. 

Now let $a$ be a smooth function defined on $T^*U$ and polynomial 
in the fibers. Then locally 
$a = \sum_{\alpha} \left( a_{x,\alpha} \circ \pi \right) \, \zeta_{x}$
with respect to a normal coordinate system at $x \in U$ and functions
$a_{x,\alpha} \in {\cal C}^{\infty} (U)$.
Define the operator $A: {\cal D} (U) \rightarrow {\cal C}^{\infty} (U)$ by
\begin{equation}
\label{DeOpA}
   f \mapsto Af \: = \: \Oppsi{a} \, (f) \: = \: \left( U \ni x \mapsto
   \ipif \int_{T_x^*X} \, a(\xi) \, \hat{\li{f}{\psi}} (\xi)
   \, d\xi \in \C \right)
\end{equation}
and check that
\begin{equation}
  A f \, (x) \, =  \, \sum_{\alpha} 
  \, a_{x,\alpha} (x) \, (-i)^{|\alpha|} \, \dera{v_x} \left[
  \psi \, (f \circ \exp{} ) \right] \, (0_x) \, = \,
  \sum_{\alpha} \, a_{x,\alpha} (x) \, (-i)^{|\alpha|} \, \dera{z_x} f \, (x)
\end{equation}
for every $x \in U$.
Therefore $Af$ is a differential operator independent of the choice of $\psi$.

However, if $a$ is an arbitrary element of the symbol space 
$\symi (U, {\rm Hom} (E,F))$, 
Eq.~(\ref{DeOpA}) still defines a continuous operator 
$A_{\psi} = \Oppsi{a} : {\cal D} (U, E) \rightarrow {\cal C}^{\infty} (U,F)$, 
which is a pseudodifferential operator but independent of the cut-off function
$\psi$ only up to smoothing operators.  Let us show this in more detail for
the scalar case, i.e.~where $E= F =X \times \C$. The general case is proven
analogously.

First choose $\mu \in \R$ such that $a \in \sym{\mu} (U)$,
and consider the kernel distribution $ {\cal D}
(U \times U) \rightarrow \C ,$ $ (f \otimes g) \mapsto 
\: \inprod{ A_{\psi} f }{g} \: = \:
\int_X \, g(x) \, A_{\psi} f(x) \, dx$ of $A_{\psi}$.
It can be written in the following way:
\begin{eqnarray}
\label{KeOsIn}
  \lefteqn{
  \inprod{g}{ A_{\psi} f} \: = } \nonumber \\
  & = & \ipifs \int_X \int_{T^*_xX} \int_{T_xX} \psi(v)
  \, g( x  ) \, f(\exp{v}) \, e^{- i \inprod{\xi}{v}} \, a(\xi) \,
  d v \, d \xi \, dx \nonumber \\
  & = & \ipifs \int_X \hspace{0.9em}\int _{T_xX} \int_{T_x^*X} \, 
  \psi (v) \, g(x) \, f(\exp{v}) \, e^{-i \inprod{\xi}{v}} \, a (\xi) \,
  d\xi \, dv \, dx;
\end{eqnarray}
where the first integral is an iterated one, the second one an oscillatory 
integral. To check that Eq.~$(\ref{KeOsIn})$ is true use a density argument 
or in other words approximate the symbol $a$ by elements $a_k \in \symmi (U)$
(in the topology of $\sym{\tilde{\mu}} (U)$ with $\tilde{\mu} > \mu$) 
and prove the statement for the $a_k$. The claim then follows from continuity 
of both sides of Eq.~(\ref{KeOsIn}) with respect to $a$.
Let $\tilde{\psi}: TX \rightarrow [0,1]$ be another cut-off function such that 
$\tilde{\psi}|_{\tilde{O}} = 1$ on a neighborhood $W' \subset W$
of the zero section of $T^*X$. We can assume $W' = \tilde{W}$ and
claim the operator $A_{\psi} - A_{\tilde{\psi}}$ to be smoothing.
For the proof it suffices to show that the oscillatory integral
\begin{equation}
\label{KeReOsIn}
  K_{A_{\psi} - A_{\tilde{\psi}}} (v) \: = \: \ipifs \int_{T_{\pi(v)}^*X} \,
  (\psi (v) - \tilde{\psi} (v)) \, e^{- i \inprod{\xi}{v}} \, 
  a(\xi) \, d \xi, \quad v \in TU,
\end{equation}
which gives an integral representation for the kernel of 
$A_{\psi} - A_{\tilde{\psi}}$ defines a smooth function
$K_{A_{\psi} - A_{\tilde{\psi}}} \in {\cal C}^{\infty} (TU)$.
The phase function $T_{\pi (v)}^*X \ni \xi \mapsto -\xi (v) \in \R$ has only 
critical points for $v=0$. Hence for $v \neq 0$ there exists a
$(-1)$-homogeneous vertical vector field $L$ on $T^*X$ such that
$L \, e^{- i \inprod{\cdot}{v}} = e^{- i \inprod{\cdot}{v}}$. 
The adjoint $L^{\dag}$
of $L$ satisfies $\left(L^{\dag} \right)^k a \in \sym{\mu - \rho k} (U)$, where
$k \in \N$. As the amplitude 
$(v,\xi) \mapsto (\psi (v) - \tilde{\psi} (v)) \, a(\xi) $
vanishes on $W' \times_U T^*U$, the equation
\begin{equation}
\label{KeReIn}
  K_{A_{\psi} - A_{\tilde{\psi}}} (v) = \ipifs \int_{T_{\pi(v)}^*X}
  \! e^{- i \inprod{\xi}{v}} \, \left( L^{\dag} \right)^k a  (\xi) \, d\xi
\end{equation}
holds for any integer $k$ fulfilling $\mu + \rho k < -{\rm dim} X$. 
Note that the integral in Eq.~(\ref{KeReIn}) unlike the one in 
Eq.~(\ref{KeReOsIn}) is to be understood in the sense of Lebesgue.
Hence $K_{A_{\psi} - A_{\tilde{\psi}}} \in {\cal C}^{\infty} (TU)$
and $A_{\psi} - A_{\tilde{\psi}}$ is smoothing.

Let us now prove that $A_{\psi}$ lies in the space 
$\pse{\mu} (U)$ of pseudodifferential operators of order
$\mu$ and type $(\rho,\delta)$ over $U$.
We have to check that $A_{\psi}$ is pseudolocal and 
with respect to some local coordinates looks 
like a pseudodifferential operator on an open set of $\R^d$.

({\it i}) $A_{\psi}$ is pseudolocal: Let $u,v \in {\cal D} (U)$
and $\supp{u} \cap \supp{v} = \emptyset$. Then the integral kernel $K$
of ${\cal C}^{\infty} (U) \ni f \mapsto v \, A_{\psi}( u f) \in 
{\cal C}^{\infty} (U)$ has the form
\begin{equation}
  K(v) \: = \: \ipifs \int_{T_{\pi (v)}^* X} \, \psi(v) \, a(\xi) \,
  v(\pi(v)) \, u(\exp{v)} \, e^{-i \inprod{\xi}{v}} \, d \xi.
\end{equation}
There exists an open neighborhood $W'\subset W$ of the zero section in $TX$ 
such that the amplitude $\psi \, (v \circ \pi) \, (u \circ \exp{}) \, a $ 
vanishes on $W' \times_U T^*X$. As for all nonzero $v \in TX$ the phase function
$\xi \mapsto - \inprod{\xi }{v}$ is noncritical, an argument similar to the one 
above for $K_{A_{\psi} - A_{\tilde{\psi}}}$ shows that $K$ is a smooth
function on $TU$, so
${\cal C}^{\infty} (U) \ni f \mapsto v \, A_{\psi}( u f) \in 
{\cal C}^{\infty} (U)$ is smoothing.

({\it ii}) 
$A_{\psi}$ is a pseudodifferential operator in appropriate coordinates:
Choose for $x \in U$ an open neighborhood $U_x \subset U$ so small
that any two points in $U_x$ can be connected by a unique geodesic.
In particular we assume that $U_x \times U_x$ lies in the range of
$(\pi , \exp{} ) : W \rightarrow X \times X$. Furthermore let us suppose
that $\psi \circ z_y (\tilde{y}) = 1$ for all $y, \tilde{y} \in U_x$.
Then we have for $f \in {\cal D}(U_x)$:
\begin{eqnarray}
\lefteqn{A_{\psi} f \, (y) \: =} \nonumber \\
   & = & \ipif \int_{T_y^* X} \, a (\xi) \, \hat{\li{f}{\psi}} (\xi)  \,
   d \xi \nonumber \\
   & = & \ipif \int_{T_x^* X} \, a \Big( 
   \underbrace{T^*_{z_x (y)}  \exp_x^{-1} }_{= \: d_x(y)} \,  (\xi) \Big) \,
   \Big| \det{ d_x (y) } \Big| \; {\cal F} (f \circ \exp{})  \big(
   d_x(y) \, (\xi) \big) \, d \xi \nonumber \\
   & = & \ipifs \int_{T_x^* X} \int_{T_yX} \, a \Big( d_x(y) \,  (\xi) \Big) \,
   \Big| \det{ d_x (y) } \Big| \, (f \circ \exp{}_y) (v) \: e^{- i 
   \inprod{d_x(y) \, (\xi)}{v}}
   \,d v \,  d \xi \nonumber \\
   & = & \ipifs \int_{T_x^* X} \int_{T_xX} \, a \Big( d_x(y) \,  (\xi) \Big) \,
   \Big| \det{ d_x (y)} \Big| \,
   \Big| \underbrace{\det{ (T \, ( z_y \circ z_x^{-1} )}}_{= \: Z_{y,x}} (v) \Big|
   \nonumber \\
\label{PsDiRe} 
   & & 
   (f \circ z_x^{-1} ) (v) \,
   e^{- i \inprod{d_x(y) \, (\xi)}{ z_y \circ z_x^{-1} (v)}}
   \,d v \,  d \xi 
\end{eqnarray}
Now the phase function $(v,\xi) \mapsto - \!
\inprod{d_x (y) (\xi)}{z_y \circ z_x^{-1} (v)}$
vanishes for $\xi \neq 0$ if and only if $v = z_x(y)$. 
Hence after possibly shrinking $U_x$ the Kuranishi trick gives a smooth function
$G_x : U_x \times T_xX \rightarrow \Iso{T_x^*X,T_y^*X}
\subset \Lin{T_x^*X,T_y^*X}$ such that
$\inprod{d_x(y) (\xi) }{ z_y \circ z_x^{-1} (v)} \: = \:
\inprod{G_x(y,v) (\xi)}{ v - z_x (y)}$
for all $(y,v) \in U_x \times T_xX$ and $\xi \in T_x^*X$.
A change of variables in
Eq.~(\ref{PsDiRe}) then implies:
\begin{eqnarray}
\lefteqn{ A_{\psi} f (y) \: = \: } 
  \nonumber \\
  & = & \ipifs  \int_{T_x^* X} \int_{T_xX} \, a \Big( \Big[
  d_x(y) \, G^{-1}_x (y,v) \Big] \,
  (\xi) \Big) \, \Big| \det{  d_x (y) } \Big| \,  \Big| Z_{y,x} (v)  \Big| \,
  \Big| \det{ G^{-1}_x(y,v) } \Big| \nonumber \\
  & & 
   (f \circ z_x^{-1} ) (v) \, e^{-i \inprod{\xi}{v - z_x(y)}} \, dv \, d\xi
  \nonumber \\
  & = & \ipifs  \int_{T_x^* X} \int_{T_xX} \, \tilde{a}_x (y,v,\xi) \,
  (f \circ z_x^{-1} ) (v) \, e^{-i \inprod{\xi}{v - z_x(y)}} \, dv \, d\xi,
\end{eqnarray}
where $(y,v,\xi) \mapsto \tilde{a}_x (y,v,\xi) =
a \Big( \Big[ d_x(y) \, G^{-1}_x(y,v) \Big] \, (\xi) \Big) \, \Big| 
\det{ d_x(y) }  \Big| \, \Big| Z_{y,x} (v) \Big|  \,  \Big| 
\det{ G^{-1}_x (y,v) }  \Big|$ 
is a symbol lying in
$\sym{\mu} (U_x \times U_x \times T_xX)$,
as $\rho + \delta \geq 1$ and $\rho > \delta$. Hence $A$ 
is a pseudodifferential operator with respect to normal coordinates 
over $U_x$ of order $\mu$ and type $(\rho,\delta)$.
\\[0.2cm]
Our considerations now lead us to
\begin{Theorem}
\label{InSyInRe}
  Let $X$ be Riemannian, $E,F$ Riemannian or Hermitian vector bundles over $X$
  and $\exp{}$ the exponential
  function corresponding to the Levi-Civita connection on $X$.
  Then any smooth cut-off function $\psi : TX \rightarrow [0,1]$
  gives rise to a linear sheaf morphism
  ${\rm Op}_{\psi} : \symi (\cdot, {\rm Hom} (E,F)) \rightarrow \psei 
  (\cdot,E,F )$, $a \mapsto A_{\psi}$ defined by
\begin{equation}
\label{DeAPsi}
  A_{\psi} f \,  (x) \: = \: \ipif \int_{T_x^*X} \, a(\xi) \,   
  \hat{\li{f}{\psi}} (\xi) \, d\xi ,
\end{equation}
  where $a \in \symi (U, {\rm Hom} (E,F))$, 
  $f \in {\cal D} (U,E)$, $ x \in U$ and $U \subset X$ open.
  This morphism preserves the natural filtrations of 
  $\symi (\cdot,{\rm Hom}(E,F))$ and $ \psei (\cdot, E,F)$. 
  In particular it maps the subsheaf $\symmi (\cdot, {\rm Hom} (E,F)) $ 
  of smoothing symbols to the subsheaf
  $ \psemi (\cdot, {\rm Hom} (E,F) ) \subset \psei (\cdot, {\rm Hom} (E,F))$
  of smoothing pseudodifferential operators.

  The quotient morphism
$ \overline{{\rm Op}} : (\symi / \symmi) (\cdot ,{\rm Hom} (E,F))
  \rightarrow (\psei  / \psemi) (\cdot ,E,F)$
  is an isomorphism and independent of the choice of the cut-off
  function $\psi$.
\end{Theorem}
\begin{Proof}
  The first part of the theorem has been shown above.
  So it remains to prove that
  $\overline{{\rm Op}}: 
  (\symi / \symmi) (U, {\rm Hom} (E,F)) \rightarrow (\psei / \psemi)(U,E,F)$
  is bijective for all open $U \subset X$. Let us postpone this
  till Theorem \ref{SyPsInRe}, where we will show that an explicit
  inverse of $\overline{{\rm Op}}$ is given by the complete symbol map 
  introduced in the following section.
\end{Proof}
By the above Theorem the effect of $\psi$ for the operator ${\rm Op}_{\psi}$ is
only a minor one, so from now on we will write ${\rm Op}$ instead of 
${\rm Op}_{\psi}$.


\section{The symbol map}
\label{SyMa}
In the sequel denote by $\varphi: X \times T^*X |_O \rightarrow \C$
the smooth phase function defined by
\begin{equation}
\label{PhFu}
X \times T^*X |_O \ni (x,\xi) \, \mapsto \: 
\inprod{\xi}{ \exp{}^{-1}_{\pi (\xi)} (x)}\: = \:
\inprod{\xi}{ z_{\pi (\xi)} (x) } \: \in \C,
\end{equation}
where $O$ is the range of $(\varrho, \exp): W \rightarrow X \times X$.
\begin{Theorem and Definition}
\label{DeCoSy}
   Let $A \in \pse{\mu} (U,E,F)$ be a pseudodifferential operator on a 
   Riemannian manifold $X$ and $\psi :T^*X \rightarrow [0,1]$ a cut-off 
   function. Then  the $\psi$-{\bf cutted symbol} of $A$ with respect to the 
   Levi-Civita connection on $X$ is the section
\begin{equation}
\label{PsCuSy}
   \sigma_{\psi} (A) = \sigma_{\psi,A}: T^*U \rightarrow \pi^* ({\rm Hom}
   (E, F)), \quad \xi \mapsto \left[ \Xi \mapsto 
   A \left( \psi_{\pi(\xi)} (\cdot)
   \, e^{ i \varphi (\cdot, \xi)} \, \tau_{(\cdot),\pi(\xi)} \Xi
   \right) \right] \, (\pi(\xi)),
\end{equation}
where for every $x \in X$ $\psi_x = \psi \circ \exp{}_x^{-1} = \psi \circ z_x$.
$\sigma_{\psi} (A)$ is an element of $\sym{\mu} (U,X)$. 
The corresponding element $\overline{\sigma} (A) $ in the quotient 
$(\sym{\mu} / \symmi) (U, {\rm Hom} (E,F))$ is called the {\bf normal symbol}
of $A$. It is independent of the choice of $\psi$.
\end{Theorem and Definition}
\begin{Proof}
  Let us first check that $\sigma_{\psi,A}$ lies in $\sym{\mu} (U, {\rm Hom}
  (E,F))$. It suffices to assume that $E$ and $F$ are trivial bundles,
  hence that $A$ is a scalar pseudodifferential operator.
  We can find a sequence $(x_{\iota})_{\iota \in \N}$ of points in $U$ and 
  coordinate patches $U_{\iota} \subset U$ with $x_{\iota} \in U_{\iota}$ 
  such that there exist normal coordinates 
  $z_{\iota} = z_{x_{\iota}} : U_{\iota} \rightarrow T_{x_{\iota}} X$ 
  on the $U_{\iota}$. We can even assume that the operator
  $A_{\iota} : {\cal D} (U_{\iota}) \rightarrow  
  {\cal C}^{\infty} (U_{\iota}),$
  $u \mapsto  A u |_{U_{\iota}} $ induces a pseudodifferential operator on 
  $T_{x_{\iota}} X$. In other words there exist symbols 
  $a_{\iota} \in \sym{\mu} (O_{\iota} \times O_{\iota} , T_{x_{\iota}}^* X )$ 
  with $O_{\iota} = z_{\iota} (U_{\iota})$
  such that $A_{\iota}$ is given by the following oscillatory integral:
\begin{equation}
  A_{\iota} u (y) \: = \: \ipifs \int \limits_{T^*_{x_{\iota}} X} \int \limits_{T_{x_{\iota}} X}
  a_{\iota} (z_{\iota}(y),v,\xi) \, (u \circ z_{\iota}^{-1}) (v) \,
  e^{-i \inprod{\xi}{ v - z_{\iota} (y)} } \, dv d\xi .
\end{equation}
  Now let  $(\phi_{\iota})$ be a partition of unity 
  subordinate to the covering $(U_{\iota})$ of $U$ and for every index $\iota$
  let $\phi_{\iota,1} , \phi_{\iota,2} \in {\cal C}^{\infty} (U)$ such that
  $\supp{\phi_{\iota,1}} \subset U_{\iota}$, $\supp{\phi_{\iota}} \cap \supp{\phi_{\iota,2}} =
  \emptyset$ and $\phi_{\iota,1} + \phi_{\iota,2} = 1$. 
  Then the  symbol $\sigma_{\psi,A}$ can be written in
  the following form:
\begin{eqnarray}
\label{ExSiLoPs}
  \sigma_{\psi,A} (\zeta) & = &
  \left[ A \Big( \psi_{\pi (\zeta)} (\cdot)  \,
  e^{i \varphi (\cdot, \zeta) } \Big) \right] (\pi (\zeta))
  \nonumber \\
  & = & \ipifs \sum \limits_{\iota} \, \left[ \phi_{\iota,1} \, A_{\iota} 
  \Big( \phi_{\iota} (\cdot ) \, \psi_{\pi (\zeta)} (\cdot) \,
  e^{i \varphi (\cdot, \zeta) } \Big) \right] (\pi (\zeta))  \nonumber \\
  & & + \: K \Big(  \, \psi_{\pi (\zeta)} (\cdot) \,
  e^{i \varphi (\cdot, \zeta) } \Big) (\pi (\zeta)),
\end{eqnarray}
  where
\begin{equation}
  K : {\cal D} (U) \rightarrow {\cal C}^{\infty} (U) \quad f \mapsto 
  \sum \limits_{\iota} \, \phi_{\iota,2} \, A_{\iota} ( \phi_{\iota} \, f), \quad
  \zeta \in T^*X_{\iota} 
\end{equation}
  is a smoothing pseudodifferential operator.
  By the  Kuranishi trick we can assume that there exists a 
  smooth function
  $G: O_{\iota} \times O_{\iota} \rightarrow 
  \Lin{T_{x_{\iota}}^* X , T_y^* X}$ 
  such that $G(v,w)$ is invertible for every $v,w \in O_{\iota}$ and
\begin{equation}
  \inprod{G\left( z_{\iota}(y),T(z_{\iota} \circ z_y^{-1}) ( v) \right) 
  (\xi)}{v }
  \: = \: \inprod{\xi}{z_{\iota} \circ z_y^{-1} ( v ) -  z_{\iota}(y)}
\end{equation}
  for $y \in U_{\iota}$, $v \in T_y X$ appropriate and 
  $\xi \in T_{x_{\iota}}^*X$.
  Several changes of variables then give the following chain of oscillatory 
  integrals for $\zeta \in T^*U_{\iota}$:
\begin{eqnarray}
\lefteqn{  
  \left[ A_{\iota} \Big( \phi_{\iota} (\cdot ) \, \psi_{\pi ( \zeta )} ( \cdot ) \,
  e^{i \varphi (\cdot , \zeta) } \Big) \right] (y) } \nonumber \\
  & = & \ipifs \int \limits_{T_{x_{\iota}}^* X} \int \limits_{T_{x_{\iota}}X}
  \phi_{\iota,1} (y) \, a_{\iota} \! \left( z_{\iota}^{-1} (y),v,\xi \right)
  \, \phi_{\iota} (z_{\iota}^{-1} (v)) \,
  \psi (z_y \circ z_x^{-1} (v) ) \nonumber \\
  & & e^{i \inprod{\zeta}{z_y \circ z_{\iota}^{-1} (v)}} \,
  e^{-i \inprod{\xi}{v - z_{\iota} (y)}} \, dv \, d\xi \nonumber \\
  & = & \ipifs \int \limits_{T_{x_{\iota}}^* X} \int \limits_{T_{y}X}
  \phi_{\iota,1} (y) \, a_{\iota} \! \left( z_{\iota}^{-1} (y),z_{\iota} \circ z_y^{-1} (v),\xi \right) \,
  \phi_{\iota} (z_y^{-1} (v)) \, \psi (v ) \nonumber \\
  & & e^{i \inprod{\zeta}{v}} \,
  e^{-i \inprod{\xi}{z_{\iota} \circ z_y^{-1} (v) - z_{\iota} (y)}} \, 
  dv \, d\xi \nonumber \\ 
  & = & \ipifs \int \limits_{T_{x_{\iota}}^* X} \int \limits_{T_{y}X}
  \phi_{\iota,1} (y) \, a_{\iota} \! \left( z_{\iota}^{-1} (y),z_{\iota} \circ z_y^{-1} (v),\xi \right) \,
  \phi_{\iota} (z_y^{-1} (v)) \, \psi (v ) \nonumber \\
  & & e^{i \inprod{\zeta}{v}} \,
  e^{-i  \inprod{G ( z_{\iota}(y),T(z_{\iota} \circ z_y^{-1}) ( v) ) (\xi)}{v}}
  \, dv \, d\xi \nonumber \\
  & = & \ipifs \int \limits_{T_{y}^* X} \int \limits_{T_{y}X}
  \phi_{\iota,1} (y) \, a_{\iota} \! \left( z_{\iota}^{-1} (y),z_{\iota} \circ z_y^{-1} (v),
  G^{-1} \left( z_{\iota}(y),T(z_{\iota} \circ z_y^{-1}) (v) \right) (\xi) \right)
  \nonumber \\
  & & \phi_{\iota} (z_y^{-1} (v)) \, \psi (v) \,
   \Big| \det{G^{-1}\left( z_{\iota}(y),T(z_{\iota} \circ z_y^{-1}) (v) \right)}  \Big| \,
  e^{i \inprod{\zeta}{v}} \, e^{-i  \inprod{\xi}{v}}
  \, dv \, d\xi \nonumber \\
  & = & \ipifs \int \limits_{T_{x_{\iota}}^* X} \int \limits_{T_{x_{\iota}}X}
  b_{\iota}(z_{\iota} (y) , v, \xi) \,
  e^{-i  \inprod{\xi - T^* (z_y \circ z_{\iota}^{-1}) (\zeta)}{v - z_{\iota}(y)}}
  \, dv \, d\xi,
\end{eqnarray}
  where the smooth function $b_{\iota}$ on 
  $O_{\iota} \times O_{\iota} \times T_{x_{\iota}}^*X$ is defined by
\begin{eqnarray}
\lefteqn{ b_{\iota} (z_{\iota} (y) ,v,\xi ) \: =} \nonumber \\
& = & \phi_{\iota,1} (y) \, 
 a_{\iota} \left( z_{\iota}^{-1} (y),z_{\iota} \circ z_y^{-1} \circ T(z_y \circ z_{\iota}^{-1}) \, (v),
 G^{-1} \Big( z_{\iota}(y),v \Big) \circ
 \left( T^*(z_{\iota} \circ z_y^{-1} \right) (\xi)  \right) \nonumber \\
 & & \phi_{\iota} (z_y^{-1} (T(z_y \circ z_{\iota}^{-1}) \, (v))) \,
 \psi (T(z_y \circ z_{\iota}^{-1}) \, (v)) \,
 \Big| \det{G^{-1} \left( z_{\iota}(y),v \right)} \Big|
\end{eqnarray}
 and lies in $\sym{\mu} (O_{\iota} \times O_{\iota} , T_{x_{\iota}}^* X)$.
 Let $B_{\iota} \in \pse{\mu} (O_{\iota}, T_{x_{\iota}}X)$ be the
 pseudodifferential operator on $O_{\iota}$ corresponding to the symbol $b_{\iota}$.
 The symbol $\sigma_{\iota}$ defined by $\sigma_{\iota} (v,\xi) =
 e^{-i \inprod{\xi}{ v}} \, B_{\iota} \left( e^{i \inprod{\xi}{\cdot}} \right)$, 
 where $v \in O_{\iota}$ and $\xi \in T_{x_{\iota}}^*X$, now is an element
 of $\sym{\mu} (O_{\iota} , T_{x_{\iota}}^* X)$ as one knows by the general
 theory of pseudodifferential operators on real vector spaces 
 (see \cite{Hor:ALPDOIII,Shu:POST,GriSjo:MADO}).
 But we have
\begin{equation}
  \left[ \phi_{\iota,1} \, A_{\iota} 
  \Big( \phi_{\iota} (\cdot ) \, \psi_{\pi (\zeta)} (\cdot) \,
  e^{i \varphi (\cdot, \zeta) } \Big) \right] (\pi (\zeta))  \: = \:
  \sigma_{\iota} \left( z_{\iota} ( \pi (\zeta )) ,
  T^* \left( z_{\pi (\zeta)} \circ z_{\iota}^{-1} \right) (\zeta) \right),
\end{equation}
  hence Eq.~(\ref{ExSiLoPs}) shows $\sigma_{\psi,A} \in \sym{\mu} (U,X)$
  if we can prove that $\sigma_{\psi ,K} \in \symmi (U,X)$
  for $\sigma_{\psi , K} (\zeta) = K 
  \Big(  \, \psi_{\pi (\zeta)} (\cdot) \,
  e^{i \varphi (\cdot, \zeta) } \Big) (\pi (\zeta)) $.
  So let $k: U \times U \rightarrow \C$ be the kernel function of $K$
  and write $\sigma_{\psi,K}$ as an integral:
\begin{equation}
\label{ExKSy}
  \sigma_{\psi , K} (\zeta ) =  \int \limits_U \, k(\pi(\zeta ),y) \,
  \psi_{\pi(\zeta )} (y) \, e^{-i \varphi (y, \zeta)} \, dy.
\end{equation}
  Let further $\Xi_1,...,\Xi_k$ be differential forms over $U$ and
  ${\cal V}_1,...,{\cal V}_k$ the vertical vector fields on $T^*U$
  corresponding to
  ${\Xi_1} ,..., { \Xi_k}$. 
  Differentiating Eq.~(\ref{ExKSy}) under the integral sign
  and using an adjointness relation then gives
\begin{eqnarray}
\lefteqn{ 
  |\zeta |^{2m} \, {\cal V}_1 \cdot ... \cdot {\cal V}_k \, \sigma_{\psi,K} \, (\zeta)
  \: = } 
  \nonumber \\
  & = & \int \limits_U \, k( \pi(\zeta ),y) \, \psi_{\pi(\zeta )} (y) \,
  \inprod{\Xi_1 }{ z_{\pi(\zeta)}(y)} \, ... \, \inprod{\Xi_k }{ z_{\pi(\zeta)}(y)}
  \nonumber \\
  & & i^{2mn+k} \frac{ \partial^{2mn} }{\partial z_{\pi(\zeta)}^{(2m,...,2m)} }
  e^{-i \inprod{\zeta}{z_{\pi (\zeta)} (\cdot) } } (y) \, dy
  \nonumber \\
  & = &
  \int \limits_U \, \left(
  \frac{ \partial^{2mn} }{\partial z_{\pi(\zeta)}^{(2m,...,2m)} } \right)^{\dagger}
  \left( k( \pi(\zeta ),\cdot) \, \psi_{\pi (\zeta )} (\cdot) \: 
  \inprod{\Xi_1 }{ z_{\pi(\zeta)}(\cdot)} \, ... \, 
  \inprod{\Xi_k }{ z_{\pi(\zeta)}(\cdot)} \right) (y) \nonumber \\
  & & i^{2mn+k}  e^{-i \inprod{\zeta}{z_{\pi (\zeta)} (y)}  }  \, dy.
\end{eqnarray}
  Therefore 
  $ |\zeta |^{2m} \, {\cal V}_1 \cdot ... \cdot {\cal V}_k \, \sigma_{\psi,K}
  \, (\zeta) $ is uniformly bounded  as long as
  $\pi(\zeta)$ varies in a compact subset of $U$.
  A similar argument for horizontal derivatives of $\sigma_{\psi,K}$
  finally proves $\sigma_{\psi,K} \in \symmi (U,X)$.
 
  Next we have to show that $\sigma_{\psi,A} - \sigma_{\tilde{\psi},A}$ is a
  smoothing symbol for any second cut-off function 
  $\tilde{\psi} : T^*X \rightarrow [0,1]$. It suffices to prove that
  $\tilde{\sigma}_{\iota}$ with
\begin{equation}
  \tilde{\sigma}_{\iota} (\zeta) = A_{\iota} \, \left(  \phi_{\iota} \, (\psi - \tilde{\psi} )
  \circ z_{\pi (\zeta) } \, e^{-i \varphi (\cdot ,\zeta) }\right)
\end{equation}
  is smoothing for every $\phi_{\iota} \in {\cal D} (U_{\iota})$. We can
  find a $(-1)$-homogeneous first order differential operator
  $L$ on the vector bundle
  $T^*_{x_{\iota}} X \times T_{x_{\iota}}X \rightarrow T_{x_{\iota}}X$
  such that for all $v \neq z_{\iota} (\pi (\zeta) )$ and $\xi \in T^*_{x_{\iota}}X$
  the equation
\begin{equation}
  L  e^{-i \inprod{\cdot}{v - z_{\iota} (\pi (\zeta))}} \, (\xi) \: = \:
  e^{-i \inprod{\xi}{v - z_{\iota} (\pi (\zeta))}}
\end{equation}
  holds (see for example {\sc Grigis, Sj\"ostrand} \cite{GriSjo:MADO} 
  Lemma 1.12 for a proof). 
  Because $(\psi - \tilde{\psi}) \circ z_{\pi(\zeta)} \circ z_{\iota}^{-1}$
  vanishes on a neighborhood of $z_{\iota} (\pi (\zeta) )$,
  the following equational chain of iterated integrals is also true:
\begin{eqnarray}
\lefteqn{ \tilde{\sigma}_{\iota} (\zeta) \: = \:} \nonumber \\
  & = &
  \ipifs \int \limits_{T_{x_{\iota}}^* X} \int \limits_{T_{x_{\iota}}X}
  a_{\iota} \! \left( z_{\iota}^{-1} (y),v,\xi \right)
  \, \phi_{\iota} (z_{\iota}^{-1} (v)) \,
  \left( \psi - \tilde{\psi} \right) (z_{\pi (\zeta)} \circ z_x^{-1} (v) ) 
  \nonumber \\
  & & e^{i \inprod{\zeta}{z_{\pi (\zeta)} \circ z_{\iota}^{-1} (v)}} \,
  e^{-i \inprod{\xi}{v - z_{\iota} (\pi (\zeta))}} \, dv \, d\xi \nonumber \\
  & = & \ipifs \int \limits_{T_{x_{\iota}}^* X} \int \limits_{T_{x_{\iota}}X}
  a_{\iota} \! \left( z_{\iota}^{-1} (y),v,\xi \right)
  \, \phi_{\iota} (z_{\iota}^{-1} (v)) \,
  \left( \psi - \tilde{\psi} \right) (z_{\pi (\zeta)} \circ z_x^{-1} (v) ) 
  \nonumber \\
  & & e^{i \inprod{\zeta}{z_{\pi (\zeta)} \circ z_{\iota}^{-1} (v)}} \,
  L^k \, e^{-i \inprod{\xi}{v - z_{\iota} (\pi (\zeta))}} \, dv \, d\xi \nonumber \\
  & = & \ipifs \int \limits_{T_{x_{\iota}}^* X} \int \limits_{T_{x_{\iota}}X} \left(
  \left( L^{\dagger} \right)^k a_{\iota} \right) \! \left( z_{\iota}^{-1} (y),v,\xi \right)
  \, \phi_{\iota} (z_{\iota}^{-1} (v)) 
  \left( \psi - \tilde{\psi} \right) (z_{\pi (\zeta)} \circ z_x^{-1} (v) ) 
  \nonumber \\
  & & e^{i \inprod{\zeta}{z_{\pi (\zeta)} \circ z_{\iota}^{-1} (v)}} \,
  e^{-i \inprod{\xi}{v - z_{\iota} (\pi (\zeta))}} \, dv \, d\xi.
\end{eqnarray}
  Note that the last  integral of this chain is to be understood in the
  sense of Lebesgue if $k \in \N$ is large enough.
  The same argument which was used for proving that the symbol 
  $\sigma_{\psi,K}$ from Eq.~(\ref{ExKSy}) is smoothing now shows 
  $\sigma_{\iota} \in \symmi (U_{\iota} , X)$.
\end{Proof}
After having defined the notion of a complete symbol, we will now give 
its essential properties in the following theorem.
\begin{Theorem}
\label{SyPsInRe} Let  $A \in \pse{\mu} (U,E,F)$ be a 
  pseudodifferential operator on a  Riemannian manifold $X$ and
  $a = \sigma_{\psi,A} \in \sym{\mu} (U, {\rm Hom} (E,F))$ be its 
  $\psi$-cutted symbol with respect to the Levi-Civita connection.
  Then $A$ and the pseudodifferential operator
  $\Oppsi{a}$ defined by Eq.~(\ref{DeAPsi}) coincide
  mo\-du\-lo smoothing operators. Moreover the sheaf morphisms
$ \overline{\sigma}: (\psei / \psemi) (\cdot,E,F) 
  \rightarrow (\symi / \symmi )(\cdot, {\rm Hom} (E,F))$
  and
$ \overline{{\rm Op}} : 
  (\symi / \symmi )(\cdot, {\rm Hom} (E,F)) \rightarrow (\psei / \psemi) 
  (\cdot,E,F)$
  are inverse to each other.
\end{Theorem}
Before proving the theorem let us first state a lemma.
\begin{Lemma}
\label{LeSmRe}
  The operator 
\begin{equation} 
  K: \, {\cal D} (U,E) \rightarrow {\cal C}^{\infty}(U,F), \quad
  f \mapsto \Bigl( U \ni x \mapsto A [ (1 - \psi_x ) f]  
  \, (x) \in F \Big) 
\end{equation} 
  is smoothing for any cut-off function $\psi : TX \rightarrow \C$.
\end{Lemma}
\begin{Proof}
  For simplicity we assume that $E=F=X\times \C$. The general case follows
  analogously. Consider the operators 
  $K_{a,\iota} :{\cal D}(U) \rightarrow {\cal C}^{\infty} (U)$
  where $K_{a,\iota} f \, (x) \: = \: 
  \phi_{\iota} (x) A [(1- \psi_a)f] (x)$ for $x \in U$, $a \in U_{\iota}$. 
  The open covering $(U_{\iota})$ of $U$ and a subordinate partition of unity 
  $(\phi_{\iota})$ are chosen such that for every point $a \in U_{\iota}$ there
  exist normal coordinates $z_a$ on $U_{\iota}$.
  We can even assume that
  $\psi (z_a (y)) = 1$ for all $a, y \in U_{\iota}$.
  Hence $\supp{ (1 - \psi_a)  } \cap \supp{\phi_{\iota}} = \emptyset$,
  and $K_{a,\iota}$ is smoothing, because $A$ is pseudolocal.
  This yields a family of smooth functions $k_{\iota} : U \times U 
  \rightarrow \C$ such that
\begin{equation}
  K_{a,\iota} f \, (x) \: = \: \int \limits_U k_{\iota} (x,y) \, (1 - \psi_a (y)) \,
  f(y) \, dy.
\end{equation}
  Now the smooth function $k \in {\cal C}^{\infty} (U \times U)$ with
  $k (x,y) \: = \: \sum_{\iota} \, k_{\iota} (x,y) \, (1 - \psi_x (y))$
  is the kernel function of $K$.
\end{Proof}
\noindent
{\sc Proof} (of Theorem \ref{SyPsInRe}):
  Check that for an arbitrary cut-off function 
  $\tilde{\psi} : TX \rightarrow  \C$
  one can find a cut-off function $\psi : TX \rightarrow  \C$ 
  such that the equation
\begin{equation}
\label{InFoTrPs}
  \psi_x (y) \, f(y) \: = \: \ipif \int \limits_{T^*_xX} \psi_x(y) \,
  e^{i \inprod{\xi}{z_x(y)}} \, \hat{\li{f}{\tilde{\psi}}} (\xi) \, 
  d \xi
\end{equation}
  holds for $x,y \in U$.
  Approximating the integral in Eq.~(\ref{InFoTrPs}) by  Riemann sums,
  we can find a sequence $(i_k)_{k \in \N}$ of natural numbers, a family 
  $(\xi_{k,\iota})_{k \in \N,1 \leq \iota \leq \iota_k} $ of elements 
  $\xi_{k,\iota} \in T^*_xX$ and a sequence 
  $(\epsilon_k)_{k \in \N}$ of positive numbers such that
\begin{eqnarray}
  & \lim \limits_{k \rightarrow \infty} \, \epsilon_k \: = \: 0, & \\
  & [ (\psi_x ) f] \, (y) \: = \: \lim \limits_{k \rightarrow \infty} \,
  \frac{\epsilon_k^n}{(2 \pi)^{n/2} } \, \sum \limits_{1 \leq i \leq i_k}
  \, \psi_x (y) \, e^{i \inprod{\xi_{k,\iota}\,}{\,z_x(y)} } \,
  \hat{\li{f}{\tilde{\psi}}} (\xi_{k,\iota}). & 
\end{eqnarray} 
  By an appropriate choice of the $\xi_{k,\iota}$ one can even assume that the
  series on the right hand side of the last equation converges as a 
  function  of $y \in U$ in the topology of ${\cal D} (U)$.
  As $A$ is continuous on ${\cal D} (U)$, we now have
\begin{eqnarray}
\lefteqn{  A [ (\psi_x ) f ] (x) \: = } \nonumber \\
  & = &  \lim \limits_{k \rightarrow \infty} \,
  \frac{\epsilon_k^n}{(2 \pi)^{n/2} } \, \sum \limits_{1 \leq i \leq i_k}
  \, A \left[ 
  \psi_x (\cdot) \, e^{i \inprod{\xi_{k,\iota}\,}{\,z_x(\cdot)} } \right]
  \! (x) \: \hat{\li{f}{\tilde{\psi}}} (\xi_{k,\iota})  \nonumber \\
  & = &  \lim \limits_{k \rightarrow \infty} \, 
  \frac{\epsilon_k^n}{(2 \pi)^{n/2} } \, \sum \limits_{1 \leq i \leq i_k} \,
  \sigma_A (\xi_{k,\iota}) \, \hat{\li{f}{\tilde{\psi}}} (\xi_{k,\iota}) \nonumber \\
  & = & \ipif \int \limits_{T_x^*X} \sigma_A (\xi) \, \hat{\li{f}{\tilde{\psi}}} 
  (\xi) \, d\xi \nonumber \\
  & = & A_{\tilde{\psi}} f.
\end{eqnarray}
On the other hand Lemma \ref{LeSmRe} provides  a smoothing operator $K$
such that for $f \in {\cal D} (U)$ and $x \in U$:
\begin{equation}
  Af \, (x) \: = \: A [ (\psi_x ) f] \, (x) + Kf \, (x) \: = \: 
  A_{\tilde{\psi}} f \, (x) + K f \, (x).
\end{equation}
Thus the first part of the theorem follows.

Now let $a \in \sym{\mu} (U, {\rm Hom} (E,F))$ be a symbol, and consider
the corresponding pseudodifferen\-tial operator 
$A = A_{\psi} \in \pse{\mu} (U,E,F)$.
After possibly altering $\psi$ we can assume that there exists a second 
cut-off function $\tilde{\psi} : T^*X \rightarrow \C$ such that 
$\tilde{\psi} |_{\supp{\psi}} = 1$. Then we have for $\zeta \in T^*_xX$
and $\Xi \in E_x$
\begin{eqnarray} 
\lefteqn{  \sigma_{A,\psi} (\zeta) \Xi \: = } \nonumber \\
 & = & \ipif \int \limits_{T^*_xX} a (\xi) \, {\cal F} \left( \li{ \left(
 \psi_x (\cdot ) e^{i \varphi (\cdot,\zeta)} \tau_{x, (\cdot)} \Xi
 \right)}{\tilde{\psi} \,} \right) (\xi) \, d \xi \nonumber \\ 
 & = & \ipif \int \limits_{T^*_xX} a (\xi) \, {\cal F} \left( \psi_x 
 (\cdot ) e^{i \varphi (\cdot,\zeta)} \Xi \right) (\xi) \, d \xi
 \nonumber \\ 
 & = & 
 \ipifs \int \limits_{T^*_xX} \int \limits_{T_xX} a (\xi) \Xi \, \psi (v) \, 
 e^{- i \inprod{\xi}{v}} \, e^{i \inprod{\zeta}{v}} \, dv \, d \xi,
\end{eqnarray}
where the last integral is an iterated one. Next choose a smooth function
$\phi : \R \rightarrow [0,1]$ such that $\supp{\phi} \subset [-2,2]$
and $\phi |_{[-1,1]} = 1$ and  define 
$a_k \in \symmi (U,{\rm Hom} (E,F))$ by
$a_k (\xi) \: = \: \phi \left( \frac{|\xi|}{k} \right)$. 
Then obviously $a_k \rightarrow a$ in
$\sym{\tilde{\mu}} (U{\rm Hom} (E,F),)$ for $\tilde{\mu} > \mu$. 
We can now write:
\begin{eqnarray}
\lefteqn{  \sigma_{A,\psi} (\zeta) \Xi \: = } \nonumber \\
  & = & \lim \limits_{k \rightarrow \infty} \ipifs  \int \limits_{T^*_xX} \int \limits_{T_xX} 
  a_k (\xi) \Xi \, \psi (v) \, e^{- i \inprod{\xi - \zeta}{v}} \, 
  dv \, d \xi  \nonumber \\
  & = & \lim \limits_{k \rightarrow \infty} a_k (\zeta) \Xi 
-
  \ipif \int \limits_{T_xX} (1 - \psi (v)) \, \left[ {\cal F}^{-1} a_k \right] \!
  (v) \, e^{i \inprod{\zeta}{v}} \, dv \nonumber \\
  & = & a( \zeta ) - \lim \limits_{k \rightarrow \infty}
  \ipif \int \limits_{T_xX} (1 - \psi (v)) \, \left[ 
  {\cal F}^{-1} \left( \left( L^{\dagger} \right)^m a_k \right) \right] \!
  (v) \, e^{i \inprod{\zeta}{v}} \, dv \nonumber \\
  & = & a(\zeta) + \lim \limits_{k \rightarrow \infty} K_m (\zeta,a_k),
\end{eqnarray} 
where $m \in \N$, $L$ is a smooth $(-1)$-homogeneous first order
differential operator on the vector bundle 
$T^*X \times_U \dot{T}U \rightarrow \dot{T}U$ such that
$L e^{i \inprod{\cdot}{v}} \, (\xi) = e^{i \inprod{\xi}{v}} $ and
\begin{equation}
K_m (\zeta,a_k) \: = \: - \ipif \int \limits_{T_xX} (1 - \psi (v)) \, \left[
{\cal F}^{-1} \left( \left( L^{\dagger} \right)^m a_k \right) \right] \! 
(-v) \, e^{i \inprod{\zeta}{v}} \, dv .
\end{equation}
Note that $\dot{T}U$ consists of all nonzero vectors $v \in TU$ and
that $ \left(  L^{\dagger} \right)^m a \in \sym{\tilde{\mu} - \rho m}$
for every $\tilde{\mu} \geq \mu$. Now fixing $\tilde{\mu} \geq \mu$ and
choosing $m$ large enough that $\rho m > \tilde{\mu} + 2 ({\rm dim} X -1)$ yields
\begin{eqnarray}
\lefteqn{  K_m (\zeta, a) \: =: \: \lim \limits_{k \rightarrow \infty}
  K_m (\zeta, a_k) \: = } \nonumber \\
  & = &  \lim \limits_{k \rightarrow \infty} 
  - \ipifs  \int \limits_{T_xX} \int \limits_{T_x^*X} 
  \left( \left( L^{\dagger}\right) ^m a_k \right) (\xi) \, (1 - \psi (v)) \,
  e^{- i \inprod{\xi - \zeta}{v}} \, d\xi \, dv   \nonumber \\
  & = &  - \ipifs  \int \limits_{T_xX} \int \limits_{T_x^*X} 
  \left( \left( L^{\dagger}\right) ^m a \right) (\xi) \, (1 - \psi (v)) \,
  e^{ i \inprod{\zeta}{v}} \, e^{- i \inprod{\xi}{v}} \, d\xi \, dv 
\end{eqnarray}
in the sense of  Lebesgue integrals.
As the phase function $\xi \mapsto - \inprod{\xi}{v}$ is noncritical
for every nonzero $v$ and $1 - \psi (v)$ vanishes in an open neighborhood
of the zero section of $TU$, a standard consideration already carried out
in preceeding proofs entails that $K_m(\cdot,a)$ is a smoothing symbol.
Henceforth $\sigma_{A,\tilde{\psi}}$ and $a$ differ by the smoothing
symbol $K_m (\cdot,a)$.
This gives the second part of the theorem.
{{\mbox{ }\hfill  $ \Box $ \parskip0.5cm \par}}
  In the following let us show that the normal symbol calculus gives rise
  to a natural transformation between the functor of pseudodifferential
  operators and the functor of symbols on the category of Riemannian
  manifolds and isometric embeddings.
  
  First we have to define these two functors. Associate to any
  Riemannian manifold $X$ the algebra ${\cal A} (X) = \psei / \psemi (X)$ of 
  pseudodifferential operators on $X$ modulo smoothing ones, and
  the space ${\cal S} (X) = \symi /\symmi (X)$ of symbols on $X$ modulo 
  smoothing ones. As $Y$ is Riemannian, we have for any isometric embedding 
  $f : X \rightarrow Y$ a natural embedding $f_* : T^*X \rightarrow T^*Y$.
  So we can define ${\cal S} (f) : {\cal S} (Y) \rightarrow {\cal S} (X)$ 
  to be the pull-back $f^*$ of smooth functions on $T^*Y$ via $f_*$.
  The morphism ${\cal A} (f) : {\cal A} (Y) \rightarrow {\cal A} (X)$
  is constructed by the following procedure. If $f$ is submersive, i.e.~a
  diffeomorphism onto its image, ${\cal A} (f)$ is just the pull-back of a pseudodifferential
  operator on $Y$ to $X$ via $f$.   
  So we now assume that $f$ is not submersive. Choose an open tubular 
  neighborhood $U$ of $f(X)$ in $Y$ and a smooth cut-off function 
  $\phi :U \rightarrow [0,1]$ with compact support and $\phi |_V = 1$ 
  on a neighborhood $V \subset U$ of $f(X)$.
  As $U$ is tubular there is a canonical projection $\Pi : U \rightarrow X$
  such that $\Pi \circ f = id_X$. Now let the pseudodifferential operator 
  $P$ represent an element of ${\cal A} (Y)$.
  We then define ${\cal A} (f) (P)$ to be the equivalence class
  of the pseudodifferential operator 
\begin{equation}
  {\cal E}' \ni u \mapsto f^* \left[ P \left( \phi ( \Pi^* u)\right) \right]
  \in {\cal E}' ,
\end{equation}
  where ${\cal E}'$ denotes distributions with compact support.
  Note that modulo smoothing operators 
  $f^* \left[ P \left( \phi ( \Pi^* u)\right) \right] $ does not depend on the
  choice of $\phi$. 
  With these definitions ${\cal A}$ and ${\cal S}$ obviously form contravariant
  functors from the category of Riemannian manifolds with isometric embeddings
  as morphisms to the category of complex vector spaces and linear mappings.
  
  The following proposition now holds.
\begin{Proposition}
  The symbol map $\sigma$ forms a natural transformation from the functor 
  ${\cal A}$ to the functor ${\cal S}$.  
\end{Proposition}
\begin{Proof}
  As $f$ is isometric, the relation
\begin{equation}
  f \circ \exp \: = \: \exp \circ Tf
\end{equation}
  is true. But then the phasefunctions on $X$ and $Y$ are related by
\begin{equation}
\label{EqPhPhIs}
  \varphi_Y ( y , f_* (\xi)) \: = \: \varphi_X ( \Pi (y), \xi),
\end{equation}
  where $y \in Y$ and $\xi \in T^*X$ with $y$ and $f (\pi (\xi))$ sufficiently
  close. Hence for a pseudodifferential operator $P$ on $Y$
\begin{eqnarray}
  \left( f^* \sigma_Y (P) \right) (\xi) & = & \left[ P \left( \psi_{f(\pi(\xi))}
  (\cdot) \, \e^{i \varphi_Y (\cdot , f_* \xi)} \right) \right] (f(\pi(\xi)))
\end{eqnarray}
  and
\begin{eqnarray}
  \left( \sigma_X f^* (P) \right) (\xi) & = & \left[ P \left( \psi_{f(\pi(\xi))}
  (\cdot) \, \phi (\cdot) \, \e^{i \varphi_X (\Pi(\cdot) ,  \xi)} \right) \right] 
  (\pi(\xi)).
\end{eqnarray}
  Eq.~(\ref{EqPhPhIs}) now implies 
\begin{equation}
  f^* \sigma_Y (P) \: = \: \sigma_X f^* (P)
\end{equation}
  which gives the claim.  
\end{Proof}


\section{The symbol of the adjoint}
In the sequel we want to derive an expression for the normal symbol 
$\sigma_{A^*}$ of the adjoint of a pseudodifferential operator $A$ in 
terms of the symbol $\sigma_A$. 

Let $E,F$ be Hermitian (Riemannian) vector bundles over the Riemannian
manifold $X$, and let $< \: , \: >_E$ resp.~$< \: , \: >$ the metric on $E$
resp.$F$. Denote by $D^E$ resp.~$D^F$ the unique torsionfree metric connection
on $E$ resp.~$F$, and by $D$ the induced metric connection on ${\rm Hom} (E,F)$
resp.~$\pi^* ({\rm Hom} (E,F))$.
Next define the function $\rho: O \rightarrow \R^+$ by
\begin{equation}
  \rho(x,y) \cdot \nu_x = \nu (y) \circ (T_v \exp_x  \times ... \times 
  T_v \exp_x) , \quad x,y \in O,  \quad v =  \exp^{-1}_x (y),
\end{equation}
where $\nu$ is the canonical volume density on $T_xX$ and $\nu_x$ its 
restriction to a volume density on $T_xX$. Now we claim that the operator 
$\Opstar{a} : {\cal D} (U,F) \rightarrow {\cal C}^{\infty} (U,E)$
defined by the iterated integral
\begin{eqnarray}
\lefteqn{\Opstar{a} g (x) \: =} \nonumber \\
  & = &  \ipifs \int_{T^*_xX} \int_{T_xX} \,
  e^{-i\inprod{\xi}{v}} \, 
  \tau^{-1}_{\exp v} \left[ a^* \left( T^*_v \exp^{-1}_x (\xi) \right)
  \, g (\exp v ) \right] {\tilde{\psi}} (v) \, 
  \rho^{-1} \! \left( \exp_x (v), x \right) \, dv \, d\xi \nonumber \\
\end{eqnarray}
is the (formal) adjoint of $\Op{\, (a)}$, 
where $a \in \sym{\mu} (U,{\rm Hom} (E,F)) $, $U \subset X$
open and $\tilde{\psi}$ is the cut-off function $TX \ni v \mapsto 
\psi \left( \exp_{\exp_x v}^{-1} (x) \right) \in [0,1]$. 
Assuming $f \in {\cal D} (V,E)$, $g \in {\cal D} (U,F)$ with $V \subset U$ 
open, $V \times V \subset O$ and $\mu$ to be sufficiently small the following 
chain of (proper) integrals holds:
\begin{eqnarray}
\label{AdSyPa}
  \lefteqn{  \int \limits_X  < f(x) ,  \Opstar{a}g \, (x) >_E 
  \, d\nu(x) \: = } \nonumber \\[2mm] \hline \nonumber \\[-2mm]
  \lefteqn{ \quad 
  \omega_x, \omega \: - \: \mbox{ canonical symplectic forms on } 
  T_xX \times T^*_xX \: \mbox{ resp. } \: T^*X}
  \nonumber \\ \hline \nonumber \\
  & = & \ipifs \int \limits_X 
  \int \limits_{T_xX \times T^*_xX} \, e^{i \inprod{\xi}{v}} \, 
  < f(x) , \tau_{\exp v} \left[ a^* \left( T^*_v \exp^{-1}_x (\xi) \right) \, 
  g (\exp v) \right] > \nonumber \\ 
  & & \hspace{3cm} \tilde{\psi} (v) \, \rho^{-1} \! 
  \left( \exp_x (v), x \right)  \, d \omega^n_x (v,\xi) d\nu (x)
  \nonumber \\[2mm] \hline \nonumber \\[-2mm]
  \lefteqn{ \quad \mbox{by } \: T^* \exp_x : T^*V \rightarrow 
  T_xX \times T^*_x X \mbox{ symplectic, } \quad x \in V}
  \nonumber \\ \hline \nonumber \\
  & = & \ipifs \int \limits_X  \int \limits_{T^*V} 
  e^{i \inprod{T^*\exp_x(\xi)}{\exp^{-1}_x (\pi(\xi))}} \,
  < \tau_{\pi(\xi),x} f (x) , a^*(\xi) \, g (\pi (\xi)) > \nonumber \\
  & & \hspace{3cm} \psi (\exp_{\pi(\xi)}^{-1} (x)) \, 
  \rho^{-1} \! \left( \pi(\xi), x \right) \, d \omega^n (\xi) \, d\nu (x) 
  \nonumber \\
  \hline \nonumber \\[-2mm]
  \lefteqn{\quad \mbox{Gau{\ss} Lemma and Fubini's Theorem}}
  \nonumber \\ \hline \nonumber \\
  & = & \ipifs \int \limits_{T^*V} \int \limits_X  
  e^{- i \inprod{\xi}{\exp^{-1}_{\pi(\xi)} (x)}} \,
  < \tau_{\pi(\xi),x} f (x) , a^*(\xi) \, g (\pi (\xi)) > \nonumber \\ 
  & & \hspace{3cm} 
  \psi (\exp_{\pi(\xi)}^{-1} (x)) \, \rho^{-1} \! \left( \pi(\xi), x \right) 
  \, d\nu (x) \, d \omega^n (\xi)
  \nonumber \\ \hline \nonumber \\[-2mm]
  \lefteqn{\quad \mbox{coordinate transformation } V \ni x \mapsto 
  w = \exp^{-1}_{\pi (\xi)} (x) \in T_{\pi (\xi)} X }
  \nonumber \\ \hline \nonumber \\
  & = & \ipifs \int \limits_{T^*V} \int \limits_{T_{\pi(\xi)}X}
  e^{- i \inprod{\xi}{w}} \, < a(\xi) \, \tau^{-1}_{\exp w} f(\exp w ),g(x)>
  \, \psi (w) \, d w \, d \omega^n (\xi)
  \nonumber \\ 
  & = &  \ipifs \int \limits_X \int \limits_{T_xX \times T^*_xX}
  e^{- i \inprod{\xi}{w}} \, < a(\xi) \, \tau^{-1}_{\exp w} f(\exp w ),g(x)>
  \, \psi (w) \, d w \, d \xi \, d \nu (x)
  \nonumber \\
  & = & \int \limits_X \, < \Op{\, (a)} f  (x) , g (x) >_F \, d\nu(x) 
\end{eqnarray} 
Recall that $\symmi (U, {\rm Hom} (E,F))$ is dense in any 
$\sym{\mu} (U, {\rm Hom} (E,F))$, $\mu \in \R$ with respect to the 
Fr\'echet topology of $\sym{\tilde{\mu}} (U, {\rm Hom} (E,F))$ 
if $\mu < \tilde{\mu}$.
By using a partition of unity and the continuity of 
$\sym{\tilde{\mu}} (U, {\rm Hom} (E,F)) \ni a \mapsto 
\Op{\, (a)} f \in {\cal C}^{\infty} (U,F)$
resp.~$\sym{\tilde{\mu}} (U, {\rm Hom} (E,F)) \ni a \mapsto \Opstar{a} g \in 
{\cal C}^{\infty} (U,E)$ for any $f \in {\cal D} (U,E)$ and 
$g \in {\cal D} (U,F)$ Eq.~(\ref{AdSyPa}) now imply that
\begin{eqnarray}
  \int \limits_X \, < f(x) , \Opstar{a} g (x) >_E  \, d\nu(x)
  & = &
  \int \limits_X \, < \Op{\, (a) } f  (x) , g (x) >_F \,d\nu(x)
\end{eqnarray}
holds for every $f \in {\cal D} (U,E)$, $g \in {\cal D} (U,F)$
and $a \in \sym{\mu} (U)$ with $\mu \in \R$. 
Thus $\Opstar{a}$ is the formal adjoint of $\Op{ \, (a)}$ indeed.

The normal symbol $\sigma_{A^*}$ of $A^* = \Opstar{a}$ now is given
by
\begin{eqnarray}
\lefteqn{\sigma_{A^*} (\xi) \, \Xi 
  \: \sim \:  \left[ A^* g ( \cdot, \xi, \Xi) \right] 
  (\pi(\xi)) \: = } \nonumber\\
  & = \ipifs \int \limits_{T_xX} \int \limits_{T^*_xX}
  e^{- i \inprod{\zeta - \xi }{v}} \, \tau^{-1}_{\exp v} \left[
  a^* \left( T^*_v \exp^{-1}_x (\zeta) \right) \, \tau_{\exp v} \Xi  \right] 
  \tilde{\psi} (v) \, \psi (v) \, \rho^{-1} \! \left( \exp_x (v), x \right) 
  \, dv \, d\zeta, \nonumber \\
\end{eqnarray}        
where $\xi \in T^*X$, $x = \pi (\xi)$, $\Xi \in F_x$ and 
$g (y, \xi, \Xi ) = \psi \left( \exp^{-1}_{\psi(\xi)}(y)
\right) \, e^{i \inprod{\xi}{\exp^{-1}_{\pi(\xi)} (y)}} \, \tau_{y,\pi (\xi)}
\Xi $ with $y \in X$. 
Thus by Theorem 3.4 of {\sc Grigis, Sj\"ostrand} \cite{GriSjo:MADO}
the following asymptotic expansion holds:
\begin{equation}
\begin{split}
  \sigma_{A^*} \, (\xi) & \sim   \sum \limits_{\alpha \in \N^d} \, 
  \frac{(-i)^{|\alpha|} }{\alpha !} \, \left. \dera{v} \right|_{v =0} 
  \left. \dera{\zeta} \right|_{\zeta = \xi} \Big(
  \left[ \tau^{-1}_{\exp v} \circ 
  a^* \left( T^*_v \exp^{-1}_x (\zeta) \right) \circ \tau_{\exp v} \right] 
  \rho^{-1} (\exp_x (v),x) \Big)
  \\
  & = \sum \limits_{\alpha \in \N^d} \, 
  \frac{(-i)^{|\alpha|} }{\alpha !} \, \left.
  \frac{D^{2 |\alpha |}}{\partial z_x^{\alpha} \,
  \partial \zeta_x^{\alpha}}  \right|_{\xi}
  \Big( a^* \, \rho^{-1} (\pi(-),x) \Big) 
  \\
  & =  \sum \limits_{\alpha, \beta \in \N^d} \, 
  \frac{(-i)^{|\alpha + \beta| } }{\alpha ! \, \beta !} 
  \left[ \left. \frac{D^{|\alpha |+|\beta|}}{
  \partial z_x^{\beta} \, \partial \zeta_x^{\alpha}} \right|_{\xi} a^*
  \right]   \left[ \left. \dera{z_x} \right|_x \rho^{-1} (\cdot, x) \right]\, ,
\end{split}
\end{equation}
where 
\begin{equation}
\begin{split}
  \left. \frac{D^{|\alpha |+|\beta|}}{
  \partial z_x^{\beta} \, \partial \zeta_x^{\alpha}} \right|_{\xi} a^*
  & =  \left. \derb{v} a \right|_{v =0} 
  \left. \dera{\zeta} a \right|_{\xi}
  \left[ \tau^{-1}_{\exp v} \circ 
   a \left( T^*_v \exp^{-1}_x (\zeta) \right) \circ \tau_{\exp v} \right],  
\end{split}
\end{equation}
are symmetrized covariant derivatives of 
$a \in {\rm S}^{\infty} (X, {\rm Hom} (E,F))$ with
respect to the apriorily chosen normal coordinates $(z_x, \zeta_x)$.

Next we will search for an expansion of $\rho$ and its derivatives.
For that first write 
\begin{equation}
\label{DeThFu}
  \frac{\partial}{\partial z_x^k} = \sum \limits_{k,l}
  \theta_k^l (x,\cdot) \, e_l, 
\end{equation}
where $(e_1,...,e_d)$ is the orthonormal frame of Eq.~(\ref{NoCoOrFr}), 
and the $\theta_k^l$ are smooth functions  on an open neighborhood of the
diagonal of $X \times X$. 
Then we have the following relation
(see for example {\sc Berline, Getzler, Vergne} \cite{BerGetVer:HKDO}, p.~36):
\begin{eqnarray}
  \theta_l^k (x,y) & = & \delta_l^k - \frac{1}{6}\sum \limits_{m,n} \, 
  {R^k }_{mln} (y) \, z_x^m (y) \, z_x^n (y) + O (|z_x (y)|^3),
\end{eqnarray}
where ${R^k }_{mln}(y)$ resp.~their ``lowered'' versions ${R_{kmln}}(y)$ 
are given by the coefficients of the curvature tensor with
respect to normal coordinates at $y \in X$, i.e.~$R \left( \der{z_y^m} \right)
= \sum \limits_{k,l,n} {R^k }_{mln}(y) \der{z_y^k} \otimes dz_y^l \otimes
dz_y^n$. But this implies
\begin{eqnarray}
  \rho (x,y) & = & \left|  \det{ \left( \theta^l_k (x,y) \right)  } \right| \:
  \nonumber \\
  & = & 1 - \frac{1}{6} \sum \limits_k \sum \limits_{m,n} {R^k }_{mkn} (x) \,
  z_x^m (y) \, z_x^n (y) + O (|z_x (y)|^3 \nonumber \\
  & = & 1 - \frac{1}{6} \sum \limits_{m,n} {\rm Ric}_{mn} (x) \, z_x^m(y) \, 
  z_x^n (y) + O \left( |z_x (y)|^3 \right).
\end{eqnarray}
Using Eq.~(\ref{SyGaLe}) $z_x^k (y) = - z_y^k (x)$ for $x$ and $y$ close enough 
we now get  
\begin{eqnarray}
  \left. \frac{\partial}{\partial z_x^k} \right|_y \rho^{-1} ( \cdot, x) & = &
  \frac{1}{3} \sum \limits_l {\rm Ric}_{kl} (x) \, z_x^l (y) +  
  O \left(|z_x (y)|^2 \right), \\[3mm]
  \left. \frac{\partial^2}{\partial z_x^k \partial z_x^l} \right|_y \rho^{-1} 
  ( \cdot, x) & = & \frac{1}{3} {\rm Ric}_{kl} (x) +  
  O \left(|z_x (y)| \right) .
\end{eqnarray}
The above results are summarized by the following theorem.
\begin{Theorem}
  The formal adjoint $A^*$ of a pseudodifferential operator 
  $A \in \pse{\mu} (U,E,F)$ between Hermitian (Riemannian)
  vector bundles $E,F$ over an open subset $U$ of a Riemannian manifold $X$ is given by
\begin{eqnarray}
\lefteqn{ A^*g \, (x)  \: = \:  \Opstar{a} g \, (x) \: =} \nonumber \\
  & = & \ipifs \int_{T^*_xX} \int_{T_xX} \,
  e^{-i\inprod{\xi}{v}} \, 
  \tau^{-1}_{\exp v} \left[ a^* \left( T^*_v \exp^{-1}_x (\xi) \right)
  \, g (\exp v ) \right] {\tilde{\psi}} (v) \, 
  \rho^{-1} \! \left( \exp_x (v), x \right) \, dv \, d\xi \nonumber \\
\end{eqnarray}
where $x \in X$, $g \in {\cal D} (U,F)$ and 
$a = \sigma_A$ is the symbol of $A$. 
The symbol $\sigma_{A^*}$  of $A^*$ has asymptotic expansion
\begin{eqnarray}
  \lefteqn{ \sigma_{A^*} (\xi) \: \sim \: 
  a^* (\xi) \:+\: \sum \limits_{\alpha \in \N^d \atop |\alpha| \geq 1}
  \left. \frac{D^{2 |\alpha|}}{\partial z_{\pi (\xi)}^{\alpha} \,
  \partial \zeta_{\pi (\xi)}^{\alpha} }
  \right|_{\xi} a^* \:-} \nonumber \\
  & & -\: \frac{1}{6} \, \sum \limits_{\alpha \in \N^d} \sum \limits_{k,l} \, 
  (-i)^{|\alpha|} \left[ \left. 
  \frac{D^{2 |\alpha|+2}}{\partial z_{\pi (\xi)}^{\alpha} \,
  \partial \zeta_{\pi (\xi)}^{\alpha}  \, \partial \zeta_{\pi(\xi)}^k  \,
  \partial \zeta_{\pi(\xi)}^l} \right|_{\xi} a^* \right] 
  \, {\rm Ric}_{kl} (\pi(\xi)) \:+ \nonumber \\
  & & +\:
  \sum \limits_{\alpha, \beta \in \N^d \atop |\beta| \geq 3} \, 
  \frac{(-i)^{|\alpha \:+\: \beta| } }{\alpha ! \, \beta !} 
  \left[ \left. 
  \frac{D^{2 |\alpha|+|\beta|}}{\partial z_{\pi (\xi)}^{\alpha} \,
  \partial \zeta_{\pi (\xi)}^{\alpha + \beta }} \right|_{\xi}  a^*  \right] 
  \left[ \left. \dera{z_{\pi(\xi)}} \right|_{\pi(\xi)} \rho^{-1} (\cdot, x) 
  \right] \nonumber \\
  & = & a^* (\xi) \:-\: i \sum \limits_k \, \left. 
  \frac{D^2}{\partial z_{\pi(\xi)}^k \,
  \partial \zeta_{\pi(\xi)}^k } \right|_{\xi}  a^* 
  \: - \: \frac{1}{2} \sum \limits_{k,l} 
  \left. \frac{D^4}{\partial z_{\pi(\xi)}^k \,
  \partial z_{\pi(\xi)}^l \, 
  \partial \zeta_{\pi(\xi)}^k \, \partial \zeta_{\pi(\xi)}^l} 
  \right|_{\xi} a^* \: - \nonumber \\
  & &  - \: \frac{1}{6}  \sum \limits_{k,l} \, \left[  \left.
  \frac{D^2}{\partial \zeta_{\pi(\xi)}^k \partial \zeta_{\pi(\xi)}^l} 
  \right|_{\xi} a^* \right] {\rm Ric}_{kl} (\pi(\xi)) \: + \:  r(\xi)
\end{eqnarray}
with $\xi \in T^*X$ and $r \in \sym{\mu - 2(\rho-\delta)}(U,{\rm Hom} (E,F))$.
\end{Theorem}


\section{Product expansions}
\label{PrEx}
Most considerations of elliptic partial differential operators as well as
applications of our normal symbol to quantization theory (see {\sc Pflaum}
\cite{Pfl:LADQ}) require an expression  for the product of two
pseudodifferential operators $A,B$ in terms of the symbols of its 
components.
In the flat case on $\R^d$ it is a well known fact 
that $\sigma_{AB}$ has an asymptotic expansion of the form
\begin{equation}
  \sigma_{AB} (x, \xi) \sim \sum \limits_{\alpha} 
  \frac{(-i)^{|\alpha|}}{\alpha !} \dera{\zeta} \sigma_A (x, \xi) \, 
  \dera{x} \sigma_B (x, \xi) .
\end{equation}
Note that strictly speaking, the product $AB$ of pseudodifferential
operators is only well-defined, if at least one of them is properly supported.
But this is only a minor set-back,
as any pseudodifferential operator is properly supported modulo a
smoothing operator and the above asymptotic expansion also describes
a symbol only up to smoothing symbols.

In the following we want to derive a similar but more complicated formula 
for the case of pseudodifferential operators on manifolds.
\begin{Proposition}
\label{SyPsFu}
  Let $A \in \pse{\mu} (U,E,F)$ be a pseudodifferential operator between
  Hermitian (Riemannian) 
  vector bundles $E,F$ over a Riemannian manifold $X$  and $a = \sigma (A)$ 
  its symbol. Then for any $f \in {\cal C}^{\infty} (U,E)$ the function
\begin{equation}
  \sigma_{A,f} : T^*U \rightarrow F, \quad 
  \xi \mapsto \left[ A \left( \psi_{\pi (\xi ) } f \,
  e^{i \varphi (\cdot, \xi)} \right) \right]  (\pi (\xi))
\end{equation} 
  is a symbol of order $\mu$ and type $(\rho,\delta)$
  and has an asymptotic expansion of the form
\begin{equation}
  \sigma_{A,f} (\xi) \sim \sum \limits_{\alpha} \, 
  \frac{1}{i^{|\alpha|} \alpha !} \, \left[ \left. 
  \ddera{\zeta_{\pi(\xi)}} \right|_{\xi} a \right] \, 
  \left[ \left. \ddera{z_{\pi (\xi)}} \right|_{\pi (\xi)} f \right]  .
\end{equation}
\end{Proposition}
\begin{Proof}
  Let $x = \pi (\xi)$, $N \in \N$ and $\psi :T^*X \rightarrow \C$ a cut-off
  function. Then the following chain of equations holds:
\begin{eqnarray} 
\lefteqn{  \sigma_{A,f} (\xi) \: = \: \left[ A \left( \psi_{\pi (\xi ) } 
  f \, e^{i \varphi (\cdot, \xi) } \right) \right] (x) } \nonumber \\
  & = & \ipif \int_{T^*_xX} a(\zeta ) \, \left[ {\cal F} \left(
  \li{f}{\psi} \, e^{i \inprod{\xi}{\cdot} } \right) \right] \! (\zeta) \,
  d\zeta \nonumber \\
  & = & \ipifs \int_{T^*_xX} \int_{T_xX} \,  e^{i \inprod{ \xi - \zeta}{v}} \,
  a(\zeta) \, \tau^{-1}_{\exp v} f(\exp{v}) \, \psi (v) \, dv \, d\zeta
  \nonumber \\ 
  & = & \ipifs \int_{T^*_xX} \int_{T_xX} \, e^{- i \inprod{ \zeta}{v}} \,
  a(\xi + \zeta) \, \tau^{-1}_{\exp v} f(\exp{v}) \, \psi (v) \, dv \, d\zeta
  \nonumber \\ 
  & = & \ipifs \sum \limits_{| \alpha | \leq N} \frac{1}{\alpha!}
  \int_{T^*_xX} \int_{T_xX} \, e^{- i \inprod{ \zeta}{v}} \ddera{\zeta_x}
  a(\xi) \, \zeta^{\alpha} \, \tau^{-1}_{\exp v} f(\exp v) \, \psi (v)
  \, dv \, d\zeta + r_N (\xi)
  \nonumber \\ 
  & = & \sum \limits_{| \alpha | \leq N}
  \frac{1}{i^{|\alpha|} \, \alpha!}
  \left[ \left. \ddera{\zeta_x} \right|_{\xi } a \right] \!  \cdot 
  \left[ \left. \ddera{z_x} \right|_{x} f \right]   + r_N (\xi).
\end{eqnarray}
  Note that in this equation all integrals are iterated ones and check
  that with Taylor's formula $r_N$ is given by
\begin{equation}
\label{DeReN}
  r_N (\xi) \: = \: \ipifs \sum \limits_{| \beta | = N+1} \frac{N+1}{\beta !}
  \int_{T^*_xX} \int_0^1 (1 -t)^N \dderb{\zeta_x}  a ( \xi +
  t \zeta ) \, dt \: \zeta^{\beta} \left[ {\cal F} \left( \li{f}{\psi} \right)
  \right] \! (\zeta) \, d\zeta.
\end{equation}
  Choosing $N$ large enough, the 
  Schwarz inequality and  Plancherel's formula now entail
\begin{eqnarray}
\label{AbReN}
\lefteqn{  || r_N (\xi) || \: \leq \: } \nonumber \\
  & \ipif \sum \limits_{|\beta | = N+1} 
  \frac{N+1}{\beta !}
  \, \left| \left| \int_0^1 (1-t)^N 
  \dderb{\zeta_x} a  ( \xi + t (\cdot )) dt \right| \right|_{L^2(T_x^*X)}
  \, \left| \left| \dderb{v_x} \left. \li{f}{\psi} \right|_{T_xX}
  \right| \right|_{L^2(T_xX)}.&
\end{eqnarray}
  But for $K \subset U$ compact and $t \in [0,1]$ there exist constants 
  $D_{N,K}, \tilde{D}_{N,K} > 0$ such that
\begin{equation}
  \left| \left| \dderb{\zeta_x}  a ( \xi + t \zeta ) \right| \right|
  \leq D_{N,K} (1 + ||\xi|| + ||\zeta|| )^{\mu - \rho (N+1)}
\end{equation}
  and
\begin{eqnarray}
  \left| \left| \int_0^1 (1 - t)^k \dderb{\zeta_x} a (\xi + t (\cdot) ) dt
  \right| \right|_{L^2(T_x^*X)} & \leq & D_{N,K} \, \Big| \Big|
  ( 1 + ||\xi|| + ||\cdot|| )^{\mu - \rho (N+1) }
  \Big| \Big|_{L^2(T_x^*X)} \nonumber \\
  &  = & \tilde{D}_{N,K} \, (1 + ||\xi|| )^{n + \mu - \rho (N+1)}
\end{eqnarray}
  for all $x \in K$, $\xi, \zeta \in T^*_x X$ and all indices
  $\beta \in \N^d$ with $| \beta | = N+1$.
  Inserting this in (\ref{AbReN})
  we can find a $C_{N,K} > 0$ such that for all $x \in K$ and 
  $\xi \in T^*_xX$ 
\begin{equation}
  \left| \left| r_N (\xi) \right| \right| \leq C_{N,K} \, 
  \Big( 1 + ||\xi|| \Big)^{n + \mu -
  \rho (N+1)} \, \sum \limits_{|\gamma| \leq N + 1} \, 
  \sup \left\{ \left| \left| \ddere{\gamma}{z_x} f (y) \right| \right| \, : \:
  y \in \supp{ \psi_{\pi(\xi)}}  \right\}.
\end{equation}
 Because of $\lim \limits_{N \rightarrow \infty} 
 {\rm dim} X + \mu + \rho (N +1) = - \infty$ and a similar consideration for the 
 derivatives of $\sigma_{A,f}$ the claim now follows.
\end{Proof}
Next let us consider two pseudodifferential operators 
$A \in \psei (U,F,G)$ and $B \in \psei (U,E,F)$ between 
vector bundles $E,F,G$ over $X$. 
Modulo smoothing operators we can assume one of them
to be properly supported, so that $AB \in \psei (U,E,G)$ exists.
Defining the extended symbol 
$b^{ext} = \sigma^{ext}_{B} : U \times T^*U \rightarrow {\rm Hom} (E,F) $
of $B$ by
\begin{equation}
  b^{ext} (x,\xi) \Xi \: = \: \psi_{\pi (\xi)} (x) \, 
  e^{- i \varphi (x,\xi)} \, \left[ B \left( \psi_{\pi (\xi)} (\cdot) \,
  e^{i \varphi (\cdot ,\xi)} \tau_{(\cdot),\pi (\xi)} \Xi \right) \right] (x),
  \quad x \in U, \: (\xi,\Xi) \in \pi^* (E|_U)
\end{equation}
the symbol $\sigma_{AB}$ of the composition $AB$ is now given 
modulo $\symmi (U,{\rm Hom} (E,G))$ by
\begin{eqnarray}
\sigma_{AB} (\xi) \Xi  & = & \left[ AB \left( \psi_{\pi (\xi)}
  e^{i \varphi (\cdot,\xi)} \tau_{(\cdot),\pi (\xi)} \Xi
  \right) \right] (\pi (\xi)) 
  \: = \nonumber \\
  & = & \left[ A \left( \psi_{\pi (\xi)} \, e^{i \varphi (\cdot,\xi)} 
  \, \sigma^{ext}_{B} (\cdot,\xi) \Xi \right) \right] (\pi (\xi)).
\end{eqnarray} 
By Proposition \ref{SyPsFu} this entails that $\sigma_{AB}$ has an 
asymptotic expansion of the form
\begin{equation} 
\label{wSyExPrPs} 
  \sigma_{AB} (\xi) \sim \sum \limits_{\alpha \in \N^n} 
  \frac{1}{i^{|\alpha|} \alpha !} \, 
  \left[ \left. \ddera{\zeta_{\pi(\xi)}} \right|_{ \xi } \sigma_A  \right] 
  \left[ \left. \ddera{z_{\pi(\xi)}} \right|_{\pi (\xi) } \sigma^{ext}_{B} 
  (\cdot,\xi) \right].
\end{equation} 
As we want to find an even more detailed expansion for $\sigma_{AB}$, 
let us calculate the derivatives 
$\left. \ddera{z_{\pi(\xi)}} \right|_{\pi(\xi)} \sigma^{ext}_B (\cdot, \xi) $ 
in the following proposition.
\begin{Proposition}
  Let $B \in \pse{\mu} (U,E,F)$ be a pseudodifferential operator
  between Hermitian (Riemannian) vector bundles over the  
  Riemannian manifold $X$ and $b = \sigma (B)$ its symbol. Then for any
  $f \in {\cal C}^{\infty} (U,E)$ the smooth function
\begin{equation}
  \sigma^{ext}_{B,f} : U \times T^*U \rightarrow F, \quad
  (x,\xi) \mapsto \psi_{\pi (\xi)} (x) \, e^{-i \varphi(x,\xi )} 
  \left[ B \left( \psi_{\pi (\xi) } f \,
  e^{i \varphi (\cdot, \xi)} \right) \right] (x)
\end{equation}
  is a symbol of order $\mu$ and type $(\rho,\delta)$ on $U \times T^*U$
  and has an asymptotic expansion of the form
\begin{eqnarray}
\label{1ClEq} 
  \sigma^{ext}_{B,f} (x,\xi) & \sim &
  \sum \limits_{k \in \N} \sum \limits_{\beta,\beta_0,...,\beta_k 
  \in \N^n \atop
  {\beta_0 + ... + \beta_k = \beta \atop |\beta_1|,...,|\beta_k| \geq 2} } 
  \, \frac{i^{k-|\beta|}}{ k! \beta_0 ! \cdot ... \cdot \beta_k!} \,
  \psi_{\pi(\xi)}(x) \left[ \left. \dderb{\zeta_{x}} 
  \right|_{d_x \varphi(\: \: ,\xi) } b \right] \nonumber \\
  & &
  \left[ \left. \ddere{\beta_0}{z_{x}} \right|_x \psi_{\pi (\xi)} f \right] 
  \, \left[ \left. \dere{\beta_1}{z_{x}} \right|_x \varphi (\cdot,\xi) 
  \right] \cdot ... \cdot
  \left[ \left. \dere{\beta_k}{z_{x}} \right|_x \varphi (\cdot,\xi) \right].
\end{eqnarray} 
Furthermore the derivative $T^*U \ni \xi \mapsto \left[ \dera{z_{\pi (\xi)}}
\sigma_{B,f}^{ext} (\cdot , \xi) \right] (\pi(\xi)) \in F$ 
with $\alpha \in \N^n$ is a symbol with values in $F$ of order 
$\mu + \delta |\alpha| $ and type $(\rho,\delta)$
on $T^*U$ and has an asymptotic expansion of the form
\begin{eqnarray}
\label{2ClEq} \lefteqn{
\left[ \ddera{z_{\pi(\xi)}} \sigma_{B,f}^{ext} (\cdot,\xi)
  \right] (\pi ( \xi )) \: \sim \:  \sum \limits_{k \in \N}
  \sum \limits_{\tilde{\alpha},\alpha_0,...,\alpha_k \in \N^n \atop
  \tilde{\alpha} + \alpha_0 + ... + \alpha_k = \alpha  }
  \sum \limits_{\beta,\beta_0,...,\beta_k \in \N^n \atop
  {\beta_0 + ... + \beta_k = \beta \atop |\beta_1|,...,|\beta_k| \geq 2} } \, 
  \frac{i^{k-|\beta|} \alpha ! }{ k! \tilde{\alpha}!
  \alpha_0 ! \cdot ... \cdot \alpha_k! \beta_0 ! \cdot ... \cdot \beta_k!}}
  \nonumber \\
  & & 
  \left\{ \left. \ddere{\tilde{\alpha}}{z_{\pi(\xi)}} \right|_{\pi(\xi)} 
  \left[ \left. \dderb{\zeta_{(-)}} \right|_{d_{(-)} \varphi (\: \:,\xi)} b
  \right] \right\} \, 
  \left\{ \left. \ddere{\alpha_0}{z_{\pi (\xi)}} \right|_{\pi(\xi)} \left[
  \left. \ddere{\beta_0}{z_{(-)}} \right|_{(-)} f \right] \right\}
  \cdot \nonumber \\
  & & 
  \left\{ \left. \dere{\alpha_1}{z_{\pi(\xi)}} \right|_{\pi(\xi)} \left[ 
  \left. \dere{\beta_1}{z_{(-)}} \right|_{(-)} \! \! \varphi (\cdot,\xi)
  \right] \right\}
  \cdot ... \cdot
  \left\{ \left. \dere{\alpha_k}{z_{\pi(\xi)}} \right|_{\pi(\xi)} \left[ 
  \left. \dere{\beta_k}{z_{(-)}} \right|_{(-)} \! \! \varphi (\cdot,\xi)
  \right]  \right\}.
  \nonumber \\ & & 
\end{eqnarray}
\end{Proposition}
\begin{rNote}
\label{NoVaDe}
  The differential operators of the form $\ddera{z_{\pi(\xi)}}$ act
  on the variables denoted by $(-)$, the differential operators
  of the form $\dderb{z_{(-)}}$ on the variables $(\cdot)$.
\end{rNote}
\begin{Proof}
  Let $\eta \in {\cal C}^{\infty} (U \times U \times T^*U)$ 
  be a smooth function such that
\begin{equation}
 \eta(x,y,\xi) \: = \: \varphi(y,\xi) - \varphi(x,\xi) -
 \inprod{d_x \varphi(\: \:,\xi)}{z_x(y)}
\end{equation}
  for $x,y \in \supp{\psi_{\pi(\xi)}}$ and such that $T^*U
  \ni \xi \mapsto \eta(x,y,\xi) \in \C$ is linear for all
  $x,y \in U$. Furthermore denote by $F \in {\cal C}^{\infty} (U \times
  U \times T^*U,E)$ the function $(x,y,\xi ) \mapsto  \psi_{\pi (\xi)}
  (y) f(y) e^{ i \eta (x,y,\xi)}$.
  By the  Leibniz rule and the definition of $\eta$ we then have
\begin{eqnarray}
  \left[ \left. \dderb{z_x} \right|_y F (x,\cdot,\xi) \right]  & = &
  \sum \limits_{0 \leq k \leq \frac{|\beta|}{2} } \sum
  \limits_{\beta_0,...,\beta_k \in \N^n \atop
  {\beta_0 + ... + \beta_k = \beta \atop |\beta_1|,...,|\beta_k| \geq 2} }
  \frac{i^{k} }{ k! \beta_0 ! \cdot ... \cdot \beta_k!}  
  \left[ \left. \ddere{\beta_0}{z_x} \right|_y \psi_{\pi (\xi)} f \right]
  \cdot \nonumber \\
  & & 
  \left[ \left. \dere{\beta_1}{z_x} \right|_y \eta (x,\cdot,\xi) \right] 
  \cdot ... \cdot \left[ \left. \dere{\beta_k}{z_x} \right|_y 
  \eta (x,\cdot,\xi) \right] \, e^{i \eta(x,y,\xi)}.
\end{eqnarray}
  As $(x,y,\xi ) \mapsto \left[ \left. \dere{\beta_j}{z_x} \right|_y \eta
  (x,\cdot,\xi) \right]  $, $0 \leq j \leq k$
  is smooth and linear with respect to
  $\xi$, the preceeding equation entails that $(x,y,\xi) \mapsto
  \left[ \left. \ddere{\beta}{z_x} \right|_y F (x,\cdot,\xi) \right] $ is 
  absolutely bounded by a function $(x,y,\xi) \mapsto C_{\beta} (x,y) \, 
  |\xi|^{\frac{|\beta|}{2}}$, where
  $C_{\beta} \in {\cal C}^{\infty} (U \times U)$ and $C_{\beta} \geq 0$.
  On the other hand we have
\begin{equation}
  \sigma_{B,f}^{ext} (x,\xi)
  \: = \: \psi_{\pi (\xi)} (x) \left[ B \left( \psi_{\pi(\xi)} (\cdot)
  f (\cdot) e^{i \eta (\cdot,x,\xi) } \,
  e^{i \inprod{d_x \varphi (\: \:,\xi)}{z_x(\cdot)}} \right) \right] (x)
\end{equation}
and by the proof of Proposition \ref{SyPsFu}
\begin{eqnarray}
\label{Star}
\lefteqn{
  \sigma_{B,f}^{ext} (x,\xi) \: = \: \psi_{\pi (\xi)} (x) \sum
  \limits_{|\beta| \leq N} \frac{i^{|\beta|}}{\beta !} 
  \left[ \left. \dderb{\xi_x} \right|_{ d_x \varphi (\: \: ,\zeta ) }b\right]\,
  \left[ \left. \dderb{z_x} \right|_x F (x,\cdot,\xi) \right]  + r_N (x,\xi)}
  \nonumber \\
  & = & \psi_{\pi (\xi)} (x)
  \sum \limits_{|\beta| \leq N}
  \sum \limits_{0 \leq k \leq \frac{|\beta|}{2}}
  \sum \limits_{\beta_0,...,\beta_k \in \N^n \atop
  {\beta_0 + ... + \beta_k = \beta \atop |\beta_1|,...,|\beta_k| \geq 2} }
  \frac{i^{k - |\beta|} }{ k! \beta_0 ! \cdot ... \cdot \beta_k!}
  \left[ \left. \dderb{\xi_x} \right|_{ d_x \varphi (\: \:,\zeta ) } b\right]
  \cdot \nonumber \\ & &
  \left[ \left. \ddere{\beta_0}{z_x} \right|_x \psi_{\pi (\xi)} f \right] \cdot
  \left[ \left. \dere{\beta_1}{z_x} \right|_x \varphi (\cdot,\xi) \right] 
  \cdot ... \cdot \left[ \left.
  \dere{\beta_k}{z_x} \right|_x \varphi (\cdot,\xi) \right] 
  + r_N (x,\xi),
\end{eqnarray}
where 
\begin{eqnarray}
\label{Triangle}
\lefteqn{
  ||r_n (x,\xi)|| } \nonumber \\
  & \leq & C_{N,K}
  (1 + ||d_x \varphi (\: \:,\xi)||)^{n + \mu - \rho (N+1)}
  \sum \limits_{| \gamma | \leq N+1} \sup \left\{ \left| \left| \left[ \left. 
  \dere{\gamma}{z_x} \right|_y  F(x,\cdot,\xi ) \right] \right| \right| \, 
  : \, y \in \supp{\psi_x} \right\} \nonumber \\
  & \leq & D_{N,K} ( 1 + ||\xi||)^{n + \mu - (\rho - 1/2 )(N + 1)}
\end{eqnarray}
  holds for compact $K \subset U$, $x \in K$, $\xi \in \left. T^*X
  \right|_K$ and constants $C_{N,K},D_{N,K} > 0$.
  Because $\rho > \frac{1}{2}$, we have $\lim_{N \rightarrow \infty}
  {\rm dim} X + \mu - (\rho - \frac{1}{2} )(N + 1) = - \infty $.
  By Eqs.~(\ref{Star}) and (\ref{Triangle}) and similar relations
  for the derivatives of $\sigma_{A,f}^{ext}$ we therefore conclude
  that Eq.~(\ref{1ClEq}) is true. The rest of the claim easily follows 
  from Leibniz rule.
\end{Proof}
The  last proposition and Eq.~(\ref{wSyExPrPs}) imply
the next theorem.
\begin{Theorem}
\label{ThPrExPsDOp}
  Let $A,B \in \psei (U,E,F) $ be two pseudodifferential operators between
  Hermitian (Riemannian) vector bundles $E,F$ over a Riemannian manifold $X$.
  Assume that one of the operators $A$ and $B$ is properly supported.
  Then the symbol $\sigma_{AB}$ of the product $AB$ has the
  following asymptotic expansion. 
\begin{eqnarray} 
\label{PrExPsDOp} 
\lefteqn{
  \sigma_{AB} \, (\xi) \: \sim \: \sum \limits_{\alpha \in \N^d} \,
  \frac{i^{- |\alpha|}}{\alpha !} 
  \left[ \left. \ddera{\zeta_{\pi(\xi)}} \right|_{\xi} \sigma_A  \right]
  \left[ \left. \ddera{z_{\pi(\xi)}} \right|_{\pi(\xi)} \sigma_B  \right]
  \: + } \nonumber \\
  & & + \: \sum \limits_{k \geq 1}
  \sum \limits_{\alpha,\tilde{\alpha},\alpha_1,...,\alpha_k \in \N^n \atop
  \tilde{\alpha} + \alpha_1 + ... + \alpha_k = \alpha }
  \sum \limits_{\beta,\beta_1,...,\beta_k \in \N^n \atop
  {\beta_1 + ... + \beta_k = \beta \atop |\beta_1|,...,|\beta_k| \geq 2} } 
  \, \frac{i^{k-|\alpha|-|\beta|} }{ k!\cdot 
  \tilde{\alpha} ! \cdot \alpha_1! \cdot ... \cdot 
  \alpha_k! \beta_1 ! \cdot ... \cdot \beta_k!}  \nonumber \\
  & & \left[ \left. \ddera{\zeta_{\pi(\xi)}} \right|_{\xi}
  \sigma_A  \right] \,
  \left\{ \left. \ddere{\tilde{\alpha}}{z_{\pi(\xi)}} \right|_{\pi(\xi)}
  \left[ \left. \dderb{\zeta_{(-)}} \right|_{ d_{(-)} \varphi (\: \: ,\xi) }
  \sigma_B   \right] \right\}  \cdot \nonumber \\
  & & \left\{ \left. \dere{\alpha_1}{z_{\pi(\xi)}} \right|_{\pi(\xi)} \left[ 
  \left. \dere{\beta_1}{z_{(-)}} \right|_{(-)} \! \! \varphi (\cdot,\xi)
  \right] \right\}
  \cdot ... \cdot
  \left\{ \left. \dere{\alpha_k}{z_{\pi(\xi)}} \right|_{\pi(\xi)} \left[ 
  \left. \dere{\beta_k}{z_{(-)}} \right|_{(-)} \! \! \varphi (\cdot,\xi)
  \right]  \right\}.
  \nonumber \\ & &
\end{eqnarray}
\end{Theorem}
A somewhat more explicit expansion up to second order is given by 
\begin{Corollary}
\label{CoPrExSy}
  If $\sigma_A$ is a symbol of order $\mu$ and $\sigma_B$ a symbol of order
  $\tilde{\mu}$ the coefficients in the asymptotic expansion of $\sigma_{AB}$
  are given up to second order by the following formula:
\begin{eqnarray}
\label{FoCoPrExSy}
\lefteqn{\sigma_{AB} (\xi) \: = } \nonumber \\
 & & 
  \sigma_A (\xi) \cdot \sigma_B (\xi) \: - \: i \sum 
  \limits_{l} \frac{D \sigma_A}{\partial \zeta_{\pi (\xi),l}} (\xi) \:
  \frac{D \sigma_B}{\partial z_{\pi(\xi)}^l } (\xi) \: - \:
  \frac{1}{2} \sum \limits_{k,l} \frac{D^2 \sigma_A}{\partial 
  \zeta_{{\pi(\xi)},k} \partial \zeta_{{\pi(\xi)},l}} (\xi) \:
  \frac{D^2 \sigma_B}{\partial 
  z_{\pi(\xi)}^k \partial z_{\pi(\xi)}^l} (\xi) \nonumber \\
  & & - \: \frac{1}{12} \sum \limits_{k,l,m,n } 
  \frac{D \sigma_A}{\partial \zeta_{{\pi(\xi)},n}} (\xi) \:
  \frac{D^2 \sigma_B}{\partial 
  \zeta_{{\pi(\xi)},l} \partial \zeta_{{\pi(\xi)},m}} (\xi) \: 
  {R^k }_{mln} (\pi(\xi)) \: \zeta_{\pi(\xi),k} (\xi) \: + \: r (\xi),
\end{eqnarray}
where $\xi \in T^*U$, $x = \pi (\xi)$, $r \in \sym{\mu + \tilde{\mu} -3} (U)$
and the $R^m{ _{nkl}}(y)$ are the coefficients of the curvature tensor with 
respect to normal coordinates at $y \in X$, i.e. $R \left( \der{z_y^m} \right)
= \sum \limits_{k,l,n} {R^k } _{mln}(y) \der{z_y^k} \otimes dz_y^l \otimes
dz_y^n$. 
\end{Corollary}
\begin{Proof}
  First define the order of a summand $s = 
  s_{\alpha,\tilde{\alpha},\alpha_1,...,\alpha_k}^{\beta,\beta_1,...,\beta_k}$
  in Eq.~(\ref{PrExPsDOp})
  by ${\rm ord} (s) = |\alpha - \tilde{\alpha} + \beta -k|$. Then it is obvious
  that $s \in \sym{\mu + \tilde{\mu} + {\rm ord} (s) \, (\rho - \delta)}
  (U)$. Therefore we have to calculate only the summands $s$ with 
  ${\rm ord} (s)  \leq 2$. To achieve this let us calculate some
  coordinate changes and their derivatives. Recall the matrix-valued function 
  $\left( \theta^k_l \right)$ of Eq.~(\ref{DeThFu}) and construct its inverse 
  $\left(\bar{\theta}^k_l \right)$. In other words  
\begin{equation}
  \sum \limits_k \theta^l_k (x,y) \, \bar{\theta}^k_m (x,y) = \delta^l_m
\end{equation}
  holds for every $l,m = 1,...,d$. We can now write
\begin{eqnarray}
  \left. \der{z^k_y} \right|_z & = & \sum \limits_{l,m} \, \theta_k^l (y,z) \,
  \bar{\theta}_l^m (x,z) \, \left. \der{z_x^m} \right|_z
\end{eqnarray}
and receive the  following formulas:
\begin{eqnarray}
 \frac{ \partial z_x^k }{\partial z_y^l} (z) & = &
 \delta^k_l \: - \: \frac{1}{6} \sum \limits_{m,n} \, {R^k }_{mln} (y)
 \: z_x^m (z) \, z^n_x (z) \: \nonumber \\ & & + \:
 \frac{1}{6} \sum \limits_{m,n} \, {R^k }_{mln} (x) z_x^m (z) \, 
 z^n_x (z) \: + \: O \left( |z_x (z) + z_y (z) |^3 \right), \\[2mm]
 \label{DoDeCoCh}
 \frac{ \partial^2 z_x^k }{\partial z_y^l \partial z_y^m} (z)
 & = &
 - \: \frac{1}{6} \sum \limits_{n} \, \left( {R^k }_{mln} (y)
 + {R^k }_{nlm} (y) \right)
 \: z^n_x (z) \nonumber \\ & & 
 + \: \frac{1}{6} \sum \limits_{n} \, \left( {R^k }_{mln} (x) 
 + {R^k }_{nlm} (x) \right) 
 z^n_x (z) \nonumber \\ & & + \: O \left( |z_x (z) + z_y (z) |^2 \right), \\[2mm]
 \label{TrDeCoCh}
 \left. \der{z_x^n} \right|_{y = z} 
 \frac{ \partial^2 z_x^k }{\partial z_y^l \partial z_y^m} (y) 
 & = &
 \frac{1}{6} \, \left( {R^k }_{mln} (x) + {R^k }_{nlm} (x) \right) 
 \: + \: O \left( |z_x (z)| \right) .
\end{eqnarray}

  Next consider the functions 
\begin{eqnarray*}
  \varphi_{\alpha \beta} : \: T^*X \rightarrow \R, & &
  \xi \mapsto \left. \dera{z_{\pi(\xi)}} \right|_{y = \pi(\xi)} 
  \frac{\partial^{|\beta|} \varphi (\cdot, \xi)}{\partial^{\beta} {z_y}^{\beta} } (y)
\end{eqnarray*}
  which appear in the summands $s$ of Eq.~(\ref{PrExPsDOp}).
  One can write the $\varphi_{\alpha\beta} $ in the form 
\begin{eqnarray}
  \varphi_{\alpha \beta} (\xi) & = & \sum \limits_k \, \zeta_{\pi(\xi),k} (\xi) \: 
  \left. \dera{z_{\pi(\xi)}} \right|_{y = \pi(\xi)}
  \frac{\partial^{|\beta|} z_{\pi(\xi)}^k}{\partial {z_y}^{\beta}} (y)
\end{eqnarray}
  Thus by Eq.~(\ref{DoDeCoCh}) 
\begin{eqnarray}
\label{PhNuBe}
  \varphi_{0\beta} (\xi) & = & 0
\end{eqnarray}
  holds for every $\beta \in \N^d$ with $|\beta| \geq 2$.
  Furthermore check that 
\begin{eqnarray}
  d_y \varphi (\cdot, \xi) & = & \sum \limits_k \zeta_{\pi(\xi),k} (\xi) \left. dz_{\pi(\xi)}^k
  \right|_y
\end{eqnarray}
  is true.

  Now we are ready to calculate the summands $s$ of order ${\rm ord} (s) \leq 2$
  by considering the following cases.
\begin{enumerate}
\item ${\rm ord }(s) = 0$: 
\begin{equation}
  s (\xi) \: = \: \sigma_A (\xi) \, \sigma_B(\xi).
\end{equation}
\item ${\rm ord} (s) = 1$, $|\alpha| = 1$, $k = 0$, $\tilde{\alpha} = \alpha$:
\begin{equation}
  s (\xi) \: = \: - \: i \sum 
  \limits_{l} \frac{D \sigma_A}{\partial \zeta_{\pi (\xi),l}} (\xi) \:
  \frac{D \sigma_B}{\partial z_{\pi(\xi)}^l } (\xi).
\end{equation}
\item ${\rm ord} (s) = 1$, $|\alpha| = 0$, $k = 1$, $|\beta|=2$:
  In this case $s =0$ because of Eq.~(\ref{PhNuBe}).
\item ${\rm ord} (s) = 2$, $|\alpha| = 2$, $k = 0$, $\tilde{\alpha} = \alpha$:
\begin{equation}
  s (\xi) \: = \: - \:
  \frac{1}{2} \sum \limits_{k,l} \frac{D^2 \sigma_A}{\partial 
  \zeta_{{\pi(\xi)},k} \partial \zeta_{{\pi(\xi)},l}} (\xi) \:
  \frac{D^2 \sigma_B}{\partial 
  z_{\pi(\xi)}^k \partial z_{\pi(\xi)}^l} (\xi) .
\end{equation}
\item ${\rm ord} (s) = 2$, $|\alpha| = 1$, $k = 1$, $|\beta| =2 $:
In this case $\tilde{\alpha} =0$ by Eq.~(\ref{PhNuBe}). Hence by 
Eq.~(\ref{TrDeCoCh})
\begin{equation}
  s (\xi) \: = \:
  - \: \frac{1}{12} \sum \limits_{k,l,m,n } 
  \frac{D \sigma_A}{\partial \zeta_{{\pi(\xi)},n}} (\xi) \:
  \frac{D^2 \sigma_B}{\partial 
  \zeta_{{\pi(\xi)},l} \partial \zeta_{{\pi(\xi)},m}} (\xi) \: 
  {R^k }_{mln} (\pi(\xi)) \: \zeta_{\pi(\xi),k} (\xi).
\end{equation}
\end{enumerate}
Summing up these summands we receive the expansion of $\sigma_{AB}$ 
up to second order. This proves the claim.
\end{Proof}
The product expansion of Theorem \ref{ThPrExPsDOp} gives rise to a 
a bilinear map $\#$ on the sheaf
$\symi / \symmi (\cdot,{\rm Hom} (F,G)) \times \symi / \symmi 
(\cdot,{\rm Hom} (E,F))$ by 
\begin{eqnarray}
  \symi (U,{\rm Hom} (F,G)) \times \symi (U,{\rm Hom} (E,F))& \rightarrow& 
  \symi / \symmi (U,{\rm Hom} (E,G)), \nonumber\\
  (a,b) & \mapsto & a \# b = \sigma_{ {\rm Op} (a) \, {\rm Op} (b) }  
\end{eqnarray}
for all $U \subset X$ open and vector bundles $E,F,G$ over $X$.
The $\#$-product is an important tool in a deformation theoretical 
approach to quantization (cf.~{\sc Pflaum} \cite{Pfl:LADQ}). 
In the sequel it will be used to study ellipticity of pseudodifferential
operators on manifolds. 

\section{Ellipticity and normal symbol calculus}
The global symbol calculus introduced in the preceding paragraphs
enables us to investigate the invertibility of pseudodifferential
operators on Riemannian manifolds.
In particular we are now able define a notion of elliptic pseudodifferential
operators which is more general than the usual notion by {\sc H\"ormander} 
\cite{Hor:ALPDOIII} or {\sc Douglis, Nirenberg} \cite{DouNir:IEESPDE}.
\begin{Definition}
\label{DeElSy}
  A symbol $a \in \symi (X, {\rm Hom} (E,F))$ on a Riemannian manifold $X$ 
  is called {\bf elliptic} if there exists 
  $b_0 \in \symi (X, {\rm Hom} (F,E)))$ such that 
\begin{equation}
\label{DeEqElSy}
  a \# b_0 - 1 \in \sym{-\varepsilon} (X) \quad \mbox{ and } \quad 
  b_0 \# a -1 \in \sym{-\varepsilon} (X) 
\end{equation}
  for an $\varepsilon > 0$.
  A symbol $a \in \sym{m} (X, {\rm Hom} (E,F))$ is called elliptic of 
  order $m$ if it is elliptic and one can find a symbol 
  $b_0 \in \sym{m} (X, {\rm Hom} (F,E)) $ fulfilling Eq.(\ref{DeEqElSy}).

  An operator $A \in \psei (X,E,F)$ is called {\bf elliptic} 
  resp.~{\bf elliptic of order} $m$, if its symbol $\sigma_A$ is elliptic 
  resp.~elliptic of order $m$.
\end{Definition}

Let us give some examples of elliptic symbols resp.~operators.
\begin{rExample}
\label{ExElSy}
\begin{enumerate}
\item
  Let $a \in {\rm S}^m_{1,0} (X, {\rm Hom} (E,F))$ be a classical symbol, 
  i.e.~assume that $a$ has an expansion of the form 
  $a \sim \sum_{j \in \N} \, a_{m-j}$ with 
  $a_{m-j} \in {\rm S}^{m-j}_{1,0} (X)$ homogeneous of order
  $m-j$. Further assume that $a$ is elliptic in the classical sense, 
  i.e.~that the principal symbol $a_m$ is invertible outside the zero section.
  By the homogeneity of $a_m$ this just means that $a_m$ is elliptic
  of order $m$. But then $a$ and ${\rm Op} (a)$ must be elliptic of order 
  $m$ as well. 
\item  
  Consider the symbols $l : T^*X \rightarrow \C$, $ \xi \mapsto ||\xi ||^2$ and
  $a : T^*X \rightarrow \C$, $ \xi \mapsto \frac{1}{1+ || \xi ||^2}$
  of Example \ref{ExSy} $(ii)$. Then the pseudodifferential operator
  corresponding to $l$ is minus the Laplacian:  $\Op (l) = - \Delta$.
  Furthermore the symbols $l$ and $1 + l $ are elliptic of order $2$, and
  the relation 
  $a \# ( 1 + l) -1 $, $ (1 + l) \# a -1 \in {\rm S}^{-1}_{1,0} (X)$
  is satisfied. 
  In case $X$ is a flat Riemannian manifold, one can calculate directly
  that modulo smoothing symbols $ a \# l = l \# a = 1$, hence
  $ \Op (a)$ is a parametrix for the differential operator 
  $\Op (1 + l) = 1 - \Delta $.
\item
  Let 
  $a_{\varphi} (\xi ) = ( 1 + ||\xi||^2)^{\varphi (\pi (\xi))}$ be the symbol
  defined in Example \ref{ExSy} $(iii)$, where $\varphi : X \rightarrow \R$
  is supposed to be smooth and bounded. Then $a$ is elliptic and
  $b (\xi ) = ( 1 + ||\xi||^2)^{- \varphi (\pi (\xi))}$ fulfills
  the relations (\ref{DeEqElSy}). In case $\varphi$ is not locally constant,
  we thus receive an example of a symbol which is not elliptic
  in the sense of {\sc H\"ormander} \cite{Hor:ALPDOIII} but elliptic
  in the sense of the above definition.
\end{enumerate}
\end{rExample}

Let us now show that for any elliptic symbol $ a \in \symi (X,{\rm Hom}
(E,F))$ one can find a symbol $b \in \symi (X,{\rm Hom} (F,E))$ such that even
\begin{equation}
\label{InElSy}
  a \# b - 1 \in \symmi (X, {\rm Hom} (F,F)) \quad \mbox{ and } 
  \quad b \# a - 1 \in \symmi (X, {\rm Hom} (E,E)).
\end{equation}
Choose $b_0$ according to Definition \ref{DeElSy} and let
$r = 1 - a \# b_0$, $ l = 1 - b_0 \# a \in \symmi (X)$. Then the symbols
$b_r = b_0 \, (1 + r + r \# r + r \# r \# r + ... )$ and 
$b_l = (1 + l + l \# l + l \# l \# l + ... ) \, b_0 $ 
are well-defined and fulfill 
$ a \# b_r = 1$ and $ b_l \# a = 1$ modulo smoothing symbols. Hence
$b_r - b_l \in \symmi (X, {\rm Hom} (F,E))$, and 
$b = b_r$ is the symbol we were looking for. 
Thus we have shown the essential part for the proof of
the following theorem.
\begin{Theorem}
\label{ThElSyInOp}
  Let $A \in \psei (X,E,F)$ be an elliptic pseudodifferential operator.
  Then there exists a parametrix $B \in \psei (X,F,E)$ for $A$, i.e.~the
  relations
  \begin{equation}
  \label{InElPsPa}
    A \, B - 1 \in \psemi (X,F,F) \quad \mbox{ and } \quad 
    B \, A - 1 \in \psemi (X,E,E)
  \end{equation}
  hold. If $A$ is elliptic of order $m$, $B$ can be chosen of order $-m$.
  On the other hand, if $A$ is an operator invertible in $\psei (X)$ modulo
  smoothing operators, then $A$ is elliptic. 
  In case $X$ is compact an elliptic operator $A \in \psei (X)$ is Fredholm.
\end{Theorem}
\begin{Proof}
  Let $a = \sigma_A$. As $a$ is elliptic one can choose $b$ such that
  (\ref{InElSy}) holds. If $a$ is elliptic of order $m$, the above
  consideration shows  that $b$ is of order $-m$.
  The operator $B = \Op (b)$ then is the parametrix for $A$. 
  If on the other hand $A = \Op (a) $ has a parametrix $B = \Op (b)$, then
  $a \# b  -1$ and $b \# a - 1 $ are smoothing, hence $A$ is elliptic.
  Now recall that a pseudodifferential operator induces continuous mappings 
  between appropriate Sobolev-completions of ${\cal C}^{\infty} (X,E)$
  resp.~${\cal C}^{\infty} (X,F)$.
  and that with respect to these Sobolev-completions any smoothing 
  pseudo\-differential operator is compact. 
  Hence the claim that for compact $X$ an elliptic $A$ is Fredholm 
  follows from (\ref{InElPsPa}).
\end{Proof}

Let us give in the following propositions a rather simple criterion
for ellipticity of order $m$ in the case of scalar symbols.
\begin{Proposition}
  A scalar symbol $a \in \sym{m} (X)$ is elliptic of order $m$, 
  if and only if for every compact set $K \subset X$ there exists 
  $C_K > 0$ such that 
\begin{equation}
  || a( \xi ) || \geq \frac{1}{C_K} ||\xi ||^m
\end{equation}
  for all $\xi \in T^*X$ with $\pi (\xi) \in K$ and $|| \xi || \geq C_K$.
\end{Proposition}
\begin{Proof}
Let us first show that the condition is sufficient.
By assumption there exists a function $b \in {\cal C}^{\infty} (T^* X)$ 
such that for every compact $K \subset X$ there is $C_K > 0 $ with
\begin{equation}
  a (\xi) \cdot b \, (\xi) = 1  \quad \mbox{ and } \quad || b (\xi) ||
  \leq C_K || \xi ||^{ -m}
\end{equation}
for $\xi \in T^*X$ with $\pi (\xi) \in K$ and $||\xi|| \geq C_K$.
Hence $a \cdot b - 1 \in \symmi (X)$. After differentiating 
the relation $r = a \cdot b -1$  in local coordinates $(x,\xi)$ of 
$T^*X$ we receive
\begin{eqnarray}
\lefteqn{
  \sup_{\xi \in T^*X |_K} \, \left| \left( 1 + || \xi ||^2 \right)^{m/2} \, 
  \frac{\partial^{|\alpha|}}{\partial x^{\alpha}} \,
  \frac{\partial^{|\beta|}}{\partial \xi^{\beta}} b (\xi) \right|
  \leq 
   C \sup_{\xi \in T^*X |_K} \, \left| a (\xi) \,
  \frac{\partial^{|\alpha|}}{\partial x^{\alpha}} \,
  \frac{\partial^{|\beta|}}{\partial \xi^{\beta}} b (\xi) \, \right| \: \leq } 
  \nonumber \\
  & \leq & C \sup_{\xi \in T^*X |_K} \, \left\{ \left| r (\xi) \right| +
  \sum_{ {\alpha_1 + \alpha_2 = \alpha \atop  \beta_1 + \beta_2 = \beta}
   \atop |\alpha_1 + \beta_1 | > 0 } 
  \left| \frac{\partial^{|\alpha_1|}}{\partial x^{\alpha_1}} \,
  \frac{\partial^{|\beta_1|}}{\partial \xi^{\beta_1}} a (\xi) \, 
  \frac{\partial^{|\alpha_2|}}{\partial x^{\alpha_2}} \,
  \frac{\partial^{|\beta_2|}}{\partial \xi^{\beta_2}} 
  b (\xi) \right| \right\},
\end{eqnarray}
hence by induction $b \in \sym{-m} (X)$ follows. By the product
expansion Eq.~(\ref{FoCoPrExSy}) this implies
\begin{equation}
  a \# b -1 = r + t_1 \quad \mbox{ and } \quad b \# a - 1 = r + t_2,
\end{equation}
where $t_{1/2} \in \sym{- \varepsilon} (X)$ with 
$\varepsilon = \min \{ (\rho -\delta) , 2 \rho -1 \} > 0$. 
Hence $a$ is an elliptic symbol.

Now we will show the converse and assume the symbol $a$ to be elliptic
of order $m$. According to Theorem \ref{ThElSyInOp} there exists 
$b \in \sym{-m} (X)$ such that ${\rm Op} (b)$ is a parametrix for 
${\rm Op} (a)$. But this implies by Eq.~(\ref{FoCoPrExSy}) that
$a b -1 \in \sym{-\varepsilon}$ for some $\varepsilon >0$,
hence $ |a (\xi) b (\xi) -1| < \frac 12$ for all $\xi \in T^*X |_K$
with $||\xi || > C_K$ , where $K \subset X$ compact and $C_K > 0$.
But then
\begin{equation}
  \frac 12 < | a (\xi) b(\xi) | < C'_K |a(\xi)| || \xi||^{-m}, \quad 
  \xi \in T^*X|_K, \: ||\xi|| > C_K,
\end{equation}
which gives the claim.
\end{Proof}

In the work of {\sc Douglis, Nirenberg} \cite{DouNir:IEESPDE} a concept of
elliptic systems of pseudodifferential operators has been introduced.
We want to show in the following that this concept fits well into our framework
of ellipticity. Let us first recall the definition by {\sc Douglis}
and {\sc Nirenberg}. Let $E = E_1 \oplus ... \oplus E_K$ and 
$F= F_1 \oplus ... \oplus F_L$ be direct sums of Riemannian (Hermitian) 
vector bundles over $X$, and $A_{kl} \in \Psi^{m_{kl}}_{1,0} (X,E_k,F_l)$
with $m_{kl} \in \R$ be classical pseudodifferential operators. 
Denote the principal symbol of $A_{kl}$ by 
$p_{kl} \in {\rm S}^{m_{kl}}_{1,0} (X,E_k,F_l) $.
{\sc Douglis, Nirenberg} \cite{DouNir:IEESPDE} now call the system
$(A_{kl})$ elliptic, if there exist homogeneous symbols 
$b_{kl} \in {\rm S}^{-m_{kl}}_{1,0} (X,E_k,F_l)$ such that
\begin{equation}
\label{DoNiEl}
  \sum_l \, p_{kl} \cdot b_{lk'} - \delta_{kk'} \in {\rm S}^{-1}_{1,0} 
  (X,E_k,F_l)
  \quad \mbox{ and } \quad 
  \sum_k \, b_{lk} \cdot p_{kl'} - \delta_{ll'} \in {\rm S}^{-1}_{1,0} 
  (X,E_k,F_l) 
\end{equation}
for all $k,k',l,l'$. 
The following proposition shows that $A = (A_{kl})$ is an elliptic 
pseudodifferential operator, hence according to Theorem \ref{ThElSyInOp} 
possesses a parametrix.
\begin{Proposition}
  Let $E = E_1 \oplus ... \oplus E_K$ and $F= F_1 \oplus ... \oplus F_L$
  be direct sums of Riemannian (Hermitian) vector bundles over $X$.
  Assume that $A = (A_{kl})$ with 
  $A_{kl} \in {\rm S}^{m_{kl}}_{1,0} (X,E,F)$ comprises
  an elliptic system of pseudodifferential operators in the sense of 
  {\sc Douglis, Nirenberg} \cite{DouNir:IEESPDE}. 
  Then the pseudodifferential operator $A$ is elliptic.
\end{Proposition}
\begin{Proof}
  Write $a_{kl} = \sigma_{A_{kl}} = p_{kl} + r_{kl}$ with
  $a_{kl} \in {\rm S}^{m_{kl}}_{1,0} (X , {\rm Hom} (E_k , F_l)$
  and 
  $r_{kl} \in {\rm S}^{m_{kl}-1}_{1,0} (X , {\rm Hom} (E_k , F_l)$.
  Furthermore let $a = (a_{kl}) \in {\rm S}^{\infty}_{1,0}(X,{\rm Hom} (E,F))$,
  $b = (b_{kl}) \in {\rm S}^{\infty}_{1,0} (X,{\rm Hom} (E,F))$ and
  $r = (r_{kl}) \in {\rm S}^{\infty}_{1,0} (X,{\rm Hom} (E,F))$. 
  By Eq.~(\ref{DoNiEl}) and the expansion Eq.~\ref{FoCoPrExSy} it now follows
  \begin{equation}
  \begin{split}
    a \# b - 1 & = 
    ( p \# b -1) + r \# b \in {\rm S}^{-1}_{1,0} (X, {\rm Hom} (E,F)),
    \\ 
    b \# a  - 1 & = 
    ( b \# p -1) + b \# r \in {\rm S}^{-1}_{1,0} (X, {\rm Hom} (E,F)).
  \end{split}
  \end{equation}
  But this proves the ellipticity of $A$.
\end{Proof}
  So with the help of the calculus of normal symbols one can directly 
  construct inverses to elliptic operators in the 
  algebra of pseudodifferential operators on a Riemannian manifold. 
  In particular after having once established a global symbol calculus
  it is possible to avoid lengthy considerations in local coordinates. 
  Thus one receives a practical and natural geometric 
  tool for handling  pseudodifferential operators on manifolds.

%
%
\newpage
\nocite{Get:POSASIT,Saf:POLC,Saf:FLBO}
\bibliographystyle{amsplain}
\addcontentsline{toc}{section}{References}
\bibliography{mrabbrev,mathlib}
\end{document}